\documentclass[a4paper,10pt]{article}
\usepackage[utf8x]{inputenc}
\usepackage{tracefnt,amsmath,tabu,array}
\usepackage{amssymb,graphicx,setspace,amsfonts,amsbsy}
\usepackage{pifont,latexsym,ifthen,amsthm,rotating,calc,textcase,booktabs,cancel,slashed}

\newtheorem{theorem}{Theorem}[section]
\newtheorem{lemma}[theorem]{Lemma}
\newtheorem{corollary}[theorem]{Corollary}
\newtheorem{proposition}[theorem]{Proposition}

\newtheorem{remark}[theorem]{Remark}
\newcommand{\filledbox}{\leavevmode
  \hbox to.77778em{%
  \hfil\vbox to.675em{\hrule width.6em height.6em}\hfil}}

\newcommand{\Rm}{{\mathbb R}}

\begin{document}
\tabulinesep=1.0mm
\title{Dynamics of 3D focusing, energy-critical wave equation with radial data}

\author{Ruipeng Shen\\
Centre for Applied Mathematics\\
Tianjin University\\
Tianjin, China
}

\maketitle

\begin{abstract}
 In this article we discuss the long-time dynamics of the radial solutions to the focusing, energy-critical wave equation in 3-dimensional space. Given a solution defined for all time $t\geq 0$, we show that the soliton resolution phenomenon happens at all times $t>0$ except for a few relatively short time intervals. The main tool is the radiation theory of wave equations and the major observation of this work is a correspondence between the energy radiation and the soliton resolution/collision behaviour of solutions. We also give a few applications of the main observation on the type II blow-up solutions and ``one pass'' theory near pure mutli-solitons. 
\end{abstract}

\section{Introduction} 

\subsection{Background}

In this work we consider the long-time behaviour of the radial solutions to the focusing, energy critical wave equation in 3-dimensional space 
\[
 \left\{\begin{array}{ll} \partial_t^2 u - \Delta u = |u|^4 u, & (x,t) \in \Rm^3 \times \Rm;  \\ (u,u_t)|_{t=0} = (u_0,u_1)\in \dot{H}^1\times L^2. & \end{array} \right. \qquad \hbox{(CP1)}
\]
For convenience we use the notation $F(u) = |u|^4 u$ in this work. The energy is conserved for all $t$ in the maximal lifespan $(-T_-,T_+)$:
\[
 E = \int_{\Rm^3} \left(\frac{1}{2}|\nabla u(x,t)|^2 + \frac{1}{2}|u_t(x,t)|^2 - \frac{1}{6}|u(x,t)|^6\right) {\rm d} x.
\]
This equation is invariant under the natural dilation. More precisely, if $u$ is a solution to (CP1), then 
\[ 
 u_\lambda = \frac{1}{\lambda^{1/2}} u\left(\frac{x}{\lambda}, \frac{t}{\lambda}\right), \qquad \lambda\in \Rm^+
\]
is also a solution to (CP1). This equation is called energy critical since the initial data of $u$ and $u_\lambda$ share the same $\dot{H}^1\times L^2$ norm.

Unlike the defocusing case $\partial_t^2 u - \Delta u = - |u|^4 u$, in which all finite-energy solutions are defined for all $t\in \Rm$ and scatter in both two time directions(see, \cite{mg1, enscatter1, enscatter2, ss1, ss2}, for instance), the long time behaviour of solutions in the focusing case are quite complicated and subtle. We give a few examples: 

\paragraph{Finite time blow-up} If the solution blows up at time $T_+ \in \Rm^+$, then we must have 
\[
 \|u\|_{L^5 L^{10} ([0,T_+)\times \Rm^3)} = + \infty. 
\]
We may further divide finite time blow-up solutions into two types:
\begin{itemize} 
 \item Type I blow-up solution satisfies
  \[
   \limsup_{t\rightarrow T_+} \|(u,u_t)\|_{\dot{H}^1\times L^2} = + \infty. 
  \]
  One typical example of type I blow-up solution can be constructed in the following way: We start by considering the following solution to (CP1) 
  \[
   u(x,t) = \left(\frac{3}{4}\right)^{1/4} (T_+ -t)^{-1/2}, 
  \]
  which blows up as $t\rightarrow T_+$. In order to construct a finite-energy solution, we apply a smooth cut-off technique and utilize the finite speed of propagation. It has been proved in Donninger \cite{stablet1} that the type I blow-up of this example is stable under a small perturbation in the energy space. 
 \item Type II blow-up solution satisfies 
 \[
  \limsup_{t\rightarrow T_+} \|(u,u_t)\|_{\dot{H}^1\times L^2} < +\infty. 
 \]
 These kinds of solutions have been constructed in Krieger-Schlag-Tataru \cite{slowblowup1}, Krieger-Schlag \cite{slowblowup2} and Donninger-Huang-Krieger-Schlag \cite{moreexamples}. The behaviour of these solutions as $t\rightarrow T_+$ will be introduced in the soliton resolution part below. The instability and the stable manifolds of the specific examples given in the first two papers above have also been discussed by Krieger \cite{stablem2}, Krieger-Nahas \cite{instable} and Burzio-Krieger \cite{stablemanit2}. Similar type II blow-up solutions in higher dimensions have been discussed in Hillairet-Rapha\"{e}l \cite{4dtypeII} and Jendrej \cite{5dtypeII}. 
\end{itemize}

\paragraph{Global solutions} In this case the solution is defined for all $t\in \Rm^+$. We give two typical types of examples. The first example is the scattering solution, i.e. there exists a linear free wave $v^+$, such that 
\[
 \lim_{t\rightarrow +\infty} \|\vec{u}(t) - \vec{v}^+(t)\|_{\dot{H}^1\times L^2} = 0. 
\]
Here $\vec{u} = (u,u_t)$ and $\vec{v}^+ = (v^+, v_t^+)$. This notation will be frequently used in this work. A combination of a fixed-point argument with suitable Strichartz estimates shows that if the initial data come with a sufficiently small $\dot{H}^1\times L^2$ norm, then the corresponding solution must be a scattering solution. Another typical example of global solution is the ground state
\[
 W(x) = \left(\frac{1}{3} + |x|^2\right)^{-1/2}. 
\]
This is a stationary solution of (CP1), i.e. a solution independent of time $t$, or a solution to the elliptic equation $-\Delta u = |u|^4 u$. In fact, all radial finite-energy stationary solutions are exactly given by 
\[
 \left\{0\right\}\cup\left\{\pm W_\lambda: \lambda \in \Rm^+\right\}. 
\]
Here $W_\lambda(x)$ is the rescaled version of $W$ defined by 
\[
 W_\lambda (x) = \frac{1}{\lambda^{1/2}} W\left(\frac{x}{\lambda}\right). 
\]
More examples of global nonscattering solutions are given in Donninger-Krieger \cite{nonscaglobal1}. 

\paragraph{Soliton resolution} Soliton resolution conjecture is one of the most important open problems in the research field of dispersive and wave equations. Soliton resolution conjecture predicts that a global solution (or type II blow-up solution) to (CP1) decomposes to a sum of decoupled solitary waves, a radiation term (a linear free wave) and a small error term, as the time tends to infinity (or the blow-up time $T_+$). In the radial case, all possible nonzero solitary waves are the ground states $\pm W_\lambda$, thus we have
\begin{equation} \label{J bubble resolution}
 \vec{u}(t) = \sum_{j=1}^J \zeta_j \left(W_{\lambda_j(t)},0\right) + \vec{v}_L(t) + o(1), \quad \lambda_1(t) \gg \lambda_2(t) \gg \cdots \gg \lambda_J(t). 
\end{equation}
Here $v_L$ is a free wave and $\zeta_j \in \{\pm 1\}$. The radial case of soliton resolution in 3-dimensional space was proved by Duyckaerts-Kenig-Merle \cite{se} by a combination of profile decomposition and channel of energy method. Duyckaerts-Kenig-Merle \cite{oddhigh}, Duyckaerts-Kenig-Martel-Merle \cite{soliton4d} and Collot-Duyckaerts-Kenig-Merle \cite{soliton6d} proved the odd higher dimensional, 4-dimensional and 6-dimensional cases respectively by following the same idea, although the argument was more complicated. Recently Jendrej-Lawrie \cite{anothersoliton} gave another proof for the radial soliton resolution in dimension $d\geq 4$ by combining a sequential soliton resolution result with a ``no return'' argument. The non-radial case of soliton resolution conjecture, however, is still an open problem, although a weaker version of it, i.e. the soliton resolution along a sequence of time, has been proved by Duyckaerts-Jia-Kenig \cite{djknonradial}. 

\paragraph{Number of Bubbles} A solution like \eqref{J bubble resolution} is usually called a $J$-bubble solution. If $v_L =0$, then we call it a pure $J$-bubble solution. The specific examples of type II blow-up solutions and global nonscattering solutions given above are all one-bubble solutions in dimension 3. Solutions with at least two bubbles have been constructed in higher dimensions. Please see, for instance, Jendrej \cite{twobubble6d}. 

\paragraph{Small energy solutions} More details about the global behaviour of solutions to (CP1) are also known if the energy $E$ is not very large. For example, Kenig-Merle \cite{kenig} introduced the compactness-rigidity argument and proved that under the assumption $E(u_0,u_1)<E(W,0)$, the solution either scatters, if $\|u_0\|_{\dot{H}^1} < \|W\|_{\dot{H}^1}$; or blows up in finite time, if $\|u_0\|_{\dot{H}^1}>\|W\|_{\dot{H}^1}$. Krieger-Wong \cite{typeIblowup} shows that these blow-up solutions are of type I. The global behaviours of solutions with the threshold energy $E(u_0,u_1) = E(W,0)$ were given in Duyckaerts-Merle \cite{threshold}. Dynamics of solutions with an energy slightly greater than the ground state were discussed in Krieger-Nakanishi-Schlag \cite{aboveground1, aboveground2}. For a probability result concerning random initial data, please refer to Kenig-Mendelson \cite{randomf}.

\subsection{Main topic and idea}

In this work we consider the 3-dimensional case with radial data. We mainly focus on global solutions defined for all time $t>0$, although the idea and some of our results apply to other situations as well. We are trying to investigate the behaviour of a solution before it reach its ``final state'' of soliton resolution, especially if it takes very long time before the ``final state''. We also give another proof of the soliton resolution conjecture in 3D radial case as a direct corollary of our main result. We start by introducing the main observation of this work. 

\paragraph{Main observation} Let $(u_0,u_1)$ be radial initial data in the form of soliton resolution
\[
 (u_0,u_1) = \sum_{j=1}^J \zeta_j (W_{\lambda_j},0) + (w_0,w_1).  
\]
Here the scale separation inequality $\lambda_{j+1}/\lambda_j \ll 1$ and the smallness condition of Strichartz norm
\[
  \|\chi_0 w_L\|_{Y(\Rm)} \doteq \left(\int_{\Rm} \left(\int_{|x|>|t|} |w_L(x,t)|^{10} {\rm d} x\right)^{1/2} {\rm d} t\right)^{1/5} \ll 1
\]
holds, where $w_L$ is the linear free wave with initial data $(w_0,w_1)$ and the notations $\chi_0$ and $Y(\Rm)$ will be discussed in details at the beginning of Section 2. Then we may apply the perturbation theory (see Lemma \ref{perturbation lemma}) on the corresponding solution of (CP1) with initial data $(u_0,u_1)$ and the approximated solution 
\[
 \sum_{j=1}^J \zeta_j W_{\lambda_j} + w_L
\]
to deduce that the solution $u$ is at least defined in the region $\{(x,t): |x|>|t|\}$ and asymptotically equivalent to a linear free wave $v_L$ with weak radiation in term of the Strichartz norm. More precisely, there exists a radial free wave $v_L$ with $\|\chi_0 v_L\|_{Y(\Rm)} \ll 1$ (as defined above), such that 
\[
 \lim_{t\rightarrow \pm \infty} \int_{|x|>|t|} |\nabla_{t,x}(u-v_L)(x,t)|^2 {\rm d} x = 0. \quad \hbox{(Asymptotical equivalence)}
\]
Here $\nabla_{t,x} = (\partial_t, \nabla)$. In fact, the free waves $v_L$ and $w_L$ are very close in the energy space. In this work we observe that roughly speaking, the inverse also holds. More precisely, given a positive integer $n$, if a radial solution $u$ to (CP1) is asymptotically equivalent to a free wave $v_L$ with sufficiently weak radiation in term of the Strichartz norm, as we described above, then $\vec{u}(0)$ can be expressed in the soliton resolution form ($J\leq n$)
\[
 \vec{u} (0) =  \sum_{j=1}^J \zeta_j (W_{\lambda_j},0) + \vec{v}_L(0) + (\tilde{w}_0,\tilde{w}_1),
\]
either in the whole space $\Rm^3$ or at least outside a ball centred at the origin. In the latter case we must have $J = n$, which means that $\vec{u}(0)$ is likely to come with a relatively large energy norm when $n$ is large. Here the error term $(\tilde{w}_0,\tilde{w}_1)$ is small in the energy space. More details can be found in Section 4.

\begin{figure}[h]
 \centering
 \includegraphics[scale=0.9]{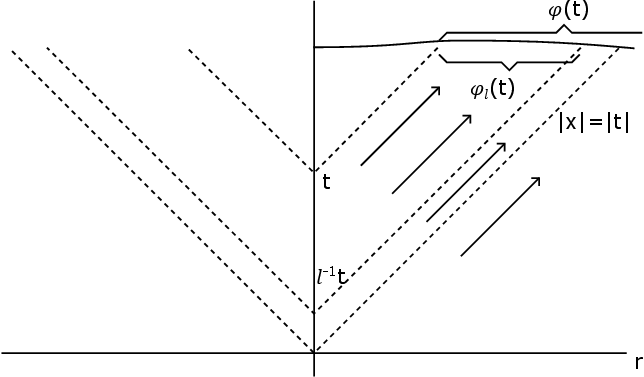}
 \caption{The radiation strength} \label{figure emission}
\end{figure}

\paragraph{General idea} Now let us describe how to use the main observation above to investigate the long time behaviour of global radial solutions to (CP1).  Let $u$ be a radial solution to (CP1) defined for all $t > 0$. It has been proved in Duyckaerts-Kenig-Merle \cite{se} that there exists a linear free wave $v_L$ with a finite energy, called the radiation part of $u$ (in the positive time direction), such that 
\begin{equation} \label{scattering part of u}
 \lim_{t\rightarrow +\infty} \int_{|x|>t - A} |\nabla_{t,x} (u-v_L)|^2 {\rm d} x = 0, \qquad \forall A \in \Rm. 
\end{equation}
The theory of radiation fields(see Section \ref{sec: radiation fields} and \ref{sec: nonlinear profiles}) implies that there exists a function $G_+ \in L^2(\Rm)$, called the radiation profile, such that the following limit holds for any fixed $A \in \Rm$:
\begin{align} \label{radiation profile u G}
 \lim_{t\rightarrow +\infty} \int_{t-A}^\infty \left(|r u_r(r,t)+G_+(r-t)|^2 + |r u_t(r,t) - G_+(r-t)|^2 \right) {\rm d} r = 0. 
\end{align}
We will try to extract information about the global behaviour of solutions from the corresponding radiation profiles $G_+$. For $t>0$, we define
\[
 \varphi (t) = \|G_+\|_{L^2([-t,+\infty))}
\]
and call it the radiation strength function. This function measures the radiation strength in the exterior region $\{(x,t'): |x|>t'-t\}$. Since a typical linear wave travels at a constant speed, it is natural to view the radiation in this region as the emission of the system before the time $t$, as illustrated in figure \ref{figure emission}. Given a large constant $\ell > 1$, we may ignore the emission before the time $\ell^{-1} t$ and focus on the emission during the time interval $[\ell^{-1}t,t]$. This gives the definition of local radiation strength function 
\[
 \varphi_{\ell} (t) = \|G_+\|_{L^2([-t, -\ell^{-1} t])}. 
\]
It follows from the explicit formula between radiation profiles and linear free wave (see Section \ref{sec: radiation fields}) that for sufficiently large constant $\ell \gg 1$ and time $t\gg 1$, the solution $u$ is asymptotically equivalent to a free wave with a weak radiation in term of the Strichartz norm in the region $\{(x,t'): |x|>|t-t'|\}$, as long as the local radiation function $\varphi_\ell (t)$ is sufficiently small at time $t$, which holds for most large time $t$ due to the fact $G_+ \in L^2(\Rm)$. An application of our main observation on the time-translated solution $u(\cdot, \cdot+t)$ then implies that the soliton resolution phenomenon happens at time $t$. Furthermore, we may verify that the soliton resolution here must hold in the whole space $\Rm^3$ by combining a continuity argument and the fact that the energy norm of $\vec{u}(t)$ can not be consistently large in a long time interval, which is a direct consequence of the virial identity. 

\subsection{Main results}

\begin{theorem} \label{main thm}
 Given any positive constants $E_0 > E(W,0)$ and $\varepsilon\ll 1$, there exists a small constant $\delta = \delta(E_0, \varepsilon)> 0$ and two large constants $\ell = \ell(E_0, \varepsilon)$, $L = L(E_0, \varepsilon)>0$ such that if $u$ is a radial solution to (CP1) satisfying
 \begin{itemize} 
  \item $u$ is defined for all time $t \geq 0$; 
  \item $R$ is a sufficiently large radius such that $\|\vec{u}(0)\|_{\dot{H}^1\times L^2(\{x: |x|>R\})} < \delta$; 
  \item the energy $E$ of $u$ satisfies $E(W,0)\leq E < E_0$;
 \end{itemize}
 then there exists a time sequence $\ell R \leq a_1 < b_1 < a_2 < b_2 < \cdots < a_m < b_m  = +\infty$ such that  
  \begin{itemize}
  \item[(a)] (Soliton resolution in stable periods) For any time interval $[a_k, b_k]$, where $k \in \{1, \cdots, m\}$, there exists a nonnegative integer $J_k$, a finite-energy free wave $v_{k,L}$, a sequence of signs $\zeta_{k,1}, \zeta_{k,2}, \cdots, \zeta_{k,J_k} \in \{\pm 1\}$ and a sequence of scale functions $\lambda_{k,1}(t) > \lambda_{k,2}(t) > \cdots > \lambda_{k,J_k}(t)$ satisfying
 \begin{align*}
  \max\left\{\frac{\lambda_{k,1}(t)}{t}, \frac{\lambda_{k,2}(t)}{\lambda_{k,1}(t)}, \cdots, \frac{\lambda_{k,J_k}(t)}{\lambda_{k,J_k-1}(t)}\right\} \leq \varepsilon^2, \qquad t\in [a_k, b_k];
 \end{align*}
 such that
  \begin{align*}
   \left\|\vec{u}(t) - \sum_{j=1}^{J_k} \zeta_{k,j} \left(W_{\lambda_{k,j}(t)},0\right) - \vec{v}_{k,L}(t)\right\|_{\dot{H}^1\times L^2} \leq \varepsilon, \qquad t\in [a_k, b_k].  
  \end{align*}
  Here $t=b_m =+\infty$ is excluded if $k=m$. In addition, the linear free wave $v_{k,L}$ satisfies 
  \begin{align*}
    &\|v_{k,L}\|_{Y([a_k,+\infty))} \leq \varepsilon; & & \|\nabla_{t,x} v_{k,L}(\cdot,t)\|_{L^2(\{x: |x|<t-\ell^{-1}a_k\})} \leq \varepsilon, \quad t\geq \ell^{-1} a_k. 
  \end{align*}  
   We call these time periods ``stable periods''. 
  \item[(b)] (Radiation concentration in collision periods) For each $k \in \{1,2,\cdots, m-1\}$, the bubble numbers $J_{k}$ and $J_{k+1}$ satisfy $J_k >J_{k+1}$. We have $\varphi_\ell (t) \geq \delta$ for each $t\in [b_{k}, a_{k+1}]$. In addition, the radiation profile $G_+$ of $u$ and times $b_k, a_{k+1}$ satisfy
  \begin{align*}
   &\left|4\pi \|G_+\|_{L^2([-a_{k+1},-b_k])}^2 - (J_k-J_{k+1}) E(W,0)\right| \leq \varepsilon^2; & & \frac{a_{k+1}}{b_k} \leq L. 
  \end{align*}
  We call these time periods ``collision periods''. In contrast, for each stable period, we have 
  \begin{align*}
   4\pi \|G_+\|_{L^2((-b_{k}, -a_k])}^2 \leq \varepsilon^2. 
  \end{align*}
 \item[(c)] (Length of preparation period) In addition, we may give a similar upper bound for the initial time of the first stable period $a_1 \leq LR$.
 \end{itemize}
\end{theorem}

\begin{remark}
 From the proof of the main theorem, we see that the radiation profile of $v_{k, L}$ in the positive time direction (see Section \ref{sec: radiation fields}) can be given by 
 \[
  G_{k,+} (s) = \left\{\begin{array}{ll} G_+ (s), & s > -b_k; \\ 0, & s < - b_k. \end{array}\right.
 \]
 In particular, the last free wave $v_{m, L}$ is exactly the scattering part $v_L$ of $u$, as defined in \eqref{scattering part of u}. It is not difficult to see that the soliton resolution conjecture is a direct consequence of Theorem \ref{main thm}. We only need to make $\varepsilon \rightarrow 0^+$. 
 \end{remark}
 
\begin{figure}[h]
 \centering
 \includegraphics[scale=1.1]{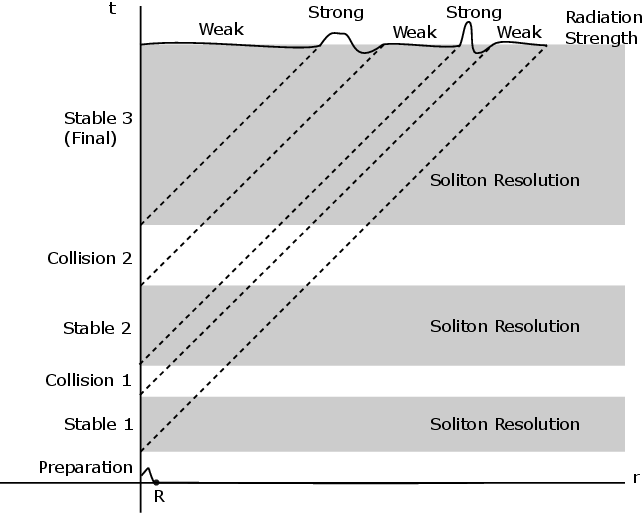}
 \caption{The relationship of stable/collision periods and radiation strength} \label{figure collision}
\end{figure}

\begin{remark}
 According to Theorem \ref{main thm}, we may split the time interval $[0,+\infty)$ into a ``preparation period'' $[0,a_1]$, several ``stable periods'' $[a_k,b_k]$, and several ``collision periods'' $[b_k, a_{k+1}]$ between consecutive stable periods. In each stable period the soliton resolution holds. It is natural to view the radiation waves travelling in the channel $t-t_2 < |x| < t - t_1$, whose strength can be measured by $\|G_+\|_{L^2([-t_2,-t_1])}^2$, as the emission of the system during the time period $[t_1,t_2]$. As a result, Theorem \ref{main thm} shows that after the preparation period, almost all radiation comes from the collision periods, whose length is bounded if we apply the logarithm transformation $t' = \ln t$. In addition, the energy of radiation in each collision period is roughly equal to the energy of bubbles eliminated in the collision. This gives a way to understand the long-time dynamics of solutions from their radiation part. This is important in physics, since the emission of energy is possibly the only thing we may actually detect for a system very far away from us. For long time dynamics of the radial 3D energy-critical wave equation, we may summarize: roughly speaking, BUBBLE COLLISION GENERATES RADIATION. Please note that a similar philosophy has been proposed in previous works, for example Duyckaerts-Kenig-Merle \cite{se}, Martel-Merle \cite{5dinelasticity} and Collot-Duyckaerts-Kenig-Merle \cite{soliton6d}, but from another perspective: briefly speaking, any solution other than ground states must send some amount of radiation, i.e. the inelasticity of bubble collision. The advantage of our method is that we may identify the time interval when the collision happens by examining the local radiation strength. Please see figure \ref{figure collision} for an illustration of stable/collision periods and their corresponding radiation strength. 
\end{remark}

\begin{remark}
 Jendrej-Lawrie \cite{anothersoliton} also introduced the conception of ``collision intervals'' in their proof, which seems similar to ``collision periods'' given the main theorem above. However, the conceptions of their ``collision intervals'' and our ``collision periods'' are actually significantly different. In Jendrej-Lawrie's proof collision intervals are time intervals in which the solution starts with but leaves the $m$-soliton state and finally comes back to the same $m$-soliton state. They eventually proved the nonexistence of collision intervals as $t\rightarrow +\infty$, which finishes their proof of soliton resolution in a combination of a sequential result of soliton resolution. In the contrast, our collision periods are time intervals in which the solution leaves from a soliton resolution state and finally reaches another different soliton resolution state (with fewer bubbles). Given an arbitrarily small $\varepsilon > 0$, collision interval exists at least for some solutions. Please see Remark \ref{existence of collision} below. 
\end{remark}

 \begin{remark} 
 For small constant $\varepsilon$, we may give an upper bound of the integer $m$: 
 \[
  m \leq J_0 + 1 \leq \left\lfloor \frac{E+\varepsilon^2}{E(W,0)} \right\rfloor + 1.
 \]
 This follows a combination of the strictly decreasing property of $J_k$ and an upper bound of $J_0$ given in the proof of the main theorem. Please see Section 5.
\end{remark}

\begin{remark} \label{existence of collision}
 Theorem \ref{main thm} is the first quantitative result of soliton resolution for energy critical wave equations, as far as the author knows. Given a global solution, we may predict the upper bound of time at which the solution first reaches an approximated soliton resolution state, simply from the energy $E$ and scale $R$ of the initial data. Of course, this soliton resolution state is not necessarily the final state of the solution. Please note that it is impossible to predict the time when the solution reaches its final state from the assumption on $u$ in Theorem \ref{main thm}. We may show this by considering a specific example. Duyckaerts-Merle \cite{threshold} constructed a radial solution to (CP1) satisfying 
 \begin{itemize}
  \item $\vec{v}(t)$ converges to $(W,0)$ in $\dot{H}^1\times L^2$ as $t\rightarrow -\infty$;
  \item $v$ scatters in the positive time direction. 
 \end{itemize}
 Thus if we choose the initial data $(u_0,u_1)$ to be $\vec{v}(t_1)$ for a large negative number $t_1$, then the initial data are closed to $(W,0)$ in $\dot{H}^1\times L^2$ but the time when the solution reaches the final(scattering) state may be arbitrarily large. This is also an example of global solutions with at least two stable periods and one collision period. Similarly Theorem \ref{main thm} also gives an upper bound on the length of each collision period between two different soliton resolution states, in term of the energy and the small constant $\varepsilon$. Please note that it is reasonable to give this upper bound by considering the quotient of two times, rather than the difference, by the scaling invariance of this equation. 
\end{remark}

\begin{remark}
 The proof of Theorem \ref{main thm} utilizes neither the profile decomposition nor a sequential version of the soliton resolution. It is much different from any previously known proof of the soliton resolution conjecture. In fact, all previously known proof heavily depends a compactness argument. The major ingredients of our proof include the radiation theory and the virial identity. The radiation theory discusses not only the energy in the exterior region, which is the main topic of the channel of energy method, but also the radiation profile defined above. The new method in this work does not involve extraction of sub-sequences of time, thus gives quantitative results (and has the potential to give more quantitative results) in a convenient way. 
\end{remark}

\begin{remark}
 Given a constant $E_0 > E(W,0)$, the dependence of constants $\delta$, $\ell$ and $L$ in the conclusion of Theorem \ref{main thm} on the small parameter $\varepsilon$ can be given explicitly
 \begin{align*}
 &\delta = c \varepsilon^2; & & \ell = E_0 c^{-2} \varepsilon^{-4}; & & L = \left(E_0 c^{-2} \varepsilon^{-4}\right)^\frac{E_0}{2\pi c^2 \varepsilon^4};  
\end{align*}
as long as $\varepsilon < \varepsilon(E_0)$ is sufficiently small. Here $c= c(E_0)$ is a constant depending on $E_0$ only. More details can be found in Section 5. 
\end{remark}

\begin{remark}
 Theorem \ref{main thm} focuses on solutions with an energy $E \geq E(W,0)$. The global behaviour of solutions with an energy $E < E(W,0)$ has been completely classified in Keng-Merle\cite{kenig}. 
 \end{remark}
 
 \begin{remark}
 Although Theorem \ref{main thm} focuses on global solutions defined for all $t\geq 0$, we may also consider evolution of solutions after a blow-up time of type II. This kind of weak solutions is called ``a canonical solution''. Please refer to Krieger-Wong \cite{canonical}. In fact, Theorem \ref{main thm} still holds for a canonical solution unless the solution eventually blows up in the manner of type I. Please see the last section of this work for more details. 
\end{remark}

\subsection{Structure of this work} 
This work is organized as follows: We introduce notations and give some preliminary results in Section 2. Then in Section 3 we show that the minimum value of $\|\vec{u}(t)\|_{\dot{H}^1\times L^2}$ in a long time interval is bounded by a constant multiple of the energy. Section 4 presents the main observation of this work: a weak radiation (in term of Strichartz norms) implies the soliton resolution phenomenon. Section 5 is devoted to the proof of the main theorem. In the final two sections we give a few more applications of our main observation given in Section 4. In Section 6 we give a ``one-pass'' theorem of pure $J$-bubble solutions; while in Section 7 we briefly discuss type II blow-up solutions and canonical solutions. 

\section{Preliminary results} 

In this section we introduce some preliminary results. We start by introducing a few useful notations for convenience. 

\paragraph{Notations} In this article we use the notation $A \lesssim B$ if there exists a constant $c$ such that the inequality $A \leq c B$ holds. Subscript(s) to the notation $\lesssim$ indicates that the constant $c$ depends on the subscript(s) but nothing else. In particular, we use $\lesssim_{1}$ to emphasize that  the constant $c$ is an absolute constant. 

\paragraph{Space norms} Assume that $R\geq 0$. We let $\mathcal{H}(R)$ be the space of restrictions of radial functions $(u_0,u_1) \in \dot{H}^1\times L^2$ to the exterior region $\{x:|x|>R\}$. The norm is given by 
\[
 \|(u_0,u_1)\|_{\mathcal{H}(R)} = \left(\int_{|x|>R} \left(|\nabla u_0(x)|^2 + |u_1(x)|^2\right){\rm d} x\right)^{1/2}.
\]
For convenience we also use the notation $\mathcal{H} = \dot{H}^1\times L^2$. 

\subsection{Exterior solutions}

For convenience of our discussion, it is helpful to introduce solutions to (CP1) defined only in an exterior region. Please note that the finite propagation speed of wave equation plays an essential role in the theory of exterior solutions. Before we discuss the basic conception of exterior solutions, we introduce a few notations. Given $R\geq 0$, we call the following region 
\[
 \Omega_R = \{(x,t)\in \Rm^3 \times \Rm: |x|>|t|+R\}
\]
an exterior region and use the notation $\chi_R$ for the characteristic function of $\Omega_R$. Similarly we let 
\[
 \Omega_{R_1,R_2} = \{(x,t)\in \Rm^3 \times \Rm: |t|+R_1<|x|<|t|+R_2\}
\]
be the channel-like region and let $\chi_{R_1,R_2}$ be its characteristic functions. Given a time interval $J$, we define the Strichartz norm 
\[
 \|u\|_{Y(J)} = \|u\|_{L^5 L^{10}(J \times \Rm^3)} = \left(\int_J \left(\int_{\Rm^3} |u(x,t)|^{10} {\rm d} x\right)^{1/2} {\rm d} t\right)^{1/5}. 
\]
These notations will be frequently used through this work. 

\paragraph{Exterior solutions} Let $u$ be a function defined in the exterior region 
\[
 \Omega =  \{(x,t): |x|>|t|+R,\, t\in (-T_1,T_2)\}.
\]
Here $T_1,T_2$ are either positive real numbers or $\infty$. We call $u$ an exterior solution to (CP1) in the region $\Omega$ with initial data $(u_0,u_1)\in \dot{H}^1 \times L^2(\Rm^3)$, if and only if $\|\chi_R u\|_{Y(J)} < +\infty$ for any bounded closed time interval $J \subset (-T_1,T_2)$ and the following identity holds:
 \[
  u = \mathbf{S}_L (u_0,u_1) + \int_0^t \frac{\sin (t-t')\sqrt{-\Delta}}{\sqrt{-\Delta}} [\chi_R(\cdot,t') F(u(\cdot,t'))] {\rm d} t', \qquad |x|>R+|t|, \; t\in (-T_1,T_2).
 \]
 The notation $\mathbf{S}_L(u_0,u_1)$ represents the linear free wave with initial data $(u_0,u_1)$. Here we multiply $u$ and $F(u)$ by the characteristic function $\chi_R$ to emphasize that $u$ and $F(u)$ are only defined in the exterior region $\Omega$. More precisely we understand the product in the following way:
 \begin{align*}
  &\chi_R u = \left\{\begin{array}{ll} u(x,t), & (x,t)\in \Omega_R; \\ 0, & (x,t)\notin \Omega_R. \end{array} \right.
  &\chi_R F(u) = \left\{\begin{array}{ll} F(u(x,t)), & (x,t) \in \Omega_R; \\ 0, & (x,t)\notin \Omega_R. \end{array} \right.
 \end{align*}
 Although we define the initial data for all $x \in \Rm^3$ in the definition above, finite speed of propagation implies that the values of initial data in the ball $\{x: |x|<R\}$ are irrelevant. As a result, when we talk about a radial exterior solution defined as above, we may specify its initial data by $(u_0,u_1)\in \mathcal{H}(R)$. 
 Similarly we may define an exterior solution $u$ to the wave equation 
 \[
  \left\{\begin{array}{l} \partial_t^2 u - \Delta u = F(x,t), \qquad (x,t)\in \Omega;\\ (u,u_t)|_{t=0} = (u_0,u_1) \end{array}\right.
 \]
 in the same manner, if $F$ is defined in the exterior region $\Omega$ and satisfies $\|\chi_R F\|_{L^1 L^2 (J \times \Rm^3)}< +\infty$ for any bounded closed interval $J \subset (-T_-,T_+)$. More details about exterior solutions can be found in Duyckaerts-Kenig-Merle\cite{oddtool}. 
 

\paragraph{Local theory} The local well-posedness of initial value problem in the exterior region immediately follows from a combination of the Strichartz estimates (see \cite{strichartz} for instance) and a fixed-point argument. The argument is similar to those in the whole space $\Rm^3$ and somewhat standard nowadays. More details of these types of argument can be found in \cite{loc1, ls}.   

\paragraph{Perturbation theory} The continuous dependence of solution on the initial data/error function immediately follows from the following lemma
\begin{lemma} \label{perturbation lemma}
 Let $M > 0$ be a constant. Then there exists two positive constants $\delta = \delta(M)$ and $C= C(M)$, such that if $v$ is a radial exterior solution to 
 \[
  \left\{\begin{array}{l} \partial_t^2 v - \Delta v = F(v) + e(x,t), \qquad (x,t)\in \Omega_R; \\
  (v,v_t)|_{t=0} = (v_0,v_1) \in \mathcal{H}(R) \end{array}\right.
 \]
 satisfying 
 \begin{align*}
  &\|\chi_R v\|_{Y(\Rm)} < M;& & \|\chi_R e(x,t)\|_{L^1 L^2(\Rm \times \Rm^3)} < \delta;
 \end{align*}
 and $(u_0,u_1)$ are a pair of radial initial data satisfying $\|(u_0,u_1)-(v_0,v_1)\|_{\mathcal{H}(R)} < \delta$, then the corresponding solution $u$ to (CP1) in the exterior region $\Omega_R$ with initial data $(u_0,u_1)$ can be defined for all $t\in \Rm$ with 
 \begin{align*}
   \|\chi_R (u-v)\|_{Y(\Rm)} + \sup_{t\in \Rm} \|\vec{u}(t) -\vec{v}(t)\|_{\mathcal{H}(R+|t|)} \leq C\left(\|\chi_R e\|_{L^1 L^2(\Rm \times \Rm^3)} + \|(u_0,u_1)-(v_0,v_1)\|_{\mathcal{H}(R)}\right). 
 \end{align*}
 Here $R\geq 0$ is an arbitrary constant. 
\end{lemma}

The proof of Lemma \ref{perturbation lemma} is the same as the situation when the solution is defined in the whole space-time $\Rm^3\times \Rm$. Please see \cite{kenig, shen2}, for instance. 

\subsection{Radiation fields of free waves} \label{sec: radiation fields} 

One of main tools of this work is the radiation field, which has a history of more than 50 years. Please see, Friedlander \cite{radiation1, radiation2} for instance. Generally speaking, radiation fields discuss the asymptotic behaviour of linear free waves. The following version of statement comes from Duyckaerts-Kenig-Merle \cite{dkm3}.

\begin{theorem}[Radiation field] \label{radiation}
Assume that $d\geq 3$ and let $u$ be a solution to the free wave equation $\partial_t^2 u - \Delta u = 0$ with initial data $(u_0,u_1) \in \dot{H}^1 \times L^2(\Rm^d)$. Then
\[
 \lim_{t\rightarrow \pm \infty} \int_{\Rm^d} \left(|\nabla u(x,t)|^2 - |u_r(x,t)|^2 + \frac{|u(x,t)|^2}{|x|^2}\right) {\rm d}x = 0
\]
 and there exist two functions $G_\pm \in L^2(\Rm \times \mathbb{S}^{d-1})$ such that
\begin{align*}
 \lim_{t\rightarrow \pm\infty} \int_0^\infty \int_{\mathbb{S}^{d-1}} \left|r^{\frac{d-1}{2}} \partial_t u(r\theta, t) - G_\pm (r\mp t, \theta)\right|^2 {\rm d}\theta {\rm d}r &= 0;\\
 \lim_{t\rightarrow \pm\infty} \int_0^\infty \int_{\mathbb{S}^{d-1}} \left|r^{\frac{d-1}{2}} \partial_r u(r\theta, t) \pm G_\pm (r\mp t, \theta)\right|^2 {\rm d}\theta {\rm d} r & = 0.
\end{align*}
In addition, the maps $(u_0,u_1) \rightarrow \sqrt{2} G_\pm$ are bijective isometries from $\dot{H}^1 \times L^2(\Rm^d)$ to $L^2 (\Rm \times \mathbb{S}^{d-1})$. 
\end{theorem}

\noindent In this work we call $G_\pm$ the radiation profiles of the linear free wave $u$, or equivalently, of the corresponding initial data $(u_0,u_1)$. Clearly the map/symmetry between radiation profiles $G_\pm$ is an isometry from $L^2(\Rm \times \mathbb{S}^{d-1})$ to itself. It is useful to give this symmetry of $G_\pm$ by an explicit formula. In this work we only need to use the 3-dimensional case (please see \cite{newradiation, shenradiation} for all dimensions, for example)
\begin{equation} \label{symmetry of G pm}
 G_+(s,\theta) = -G_-(-s, -\theta). 
\end{equation}
Thus we may uniquely determine a linear free wave by its radiation profile $G_+ \in L^2(\Rm\times \mathbb S^2)$, or radiation profile $G_- \in L^2(\Rm \times \mathbb S^2)$, or both two radiation profiles $G_\pm(s,\theta) \in L^2(\Rm^+ \times \mathbb S^2)$ (for positive $s$ only). It is not difficult to see that the free wave is radial if and only if its radiation profiles are independent of the angle $\theta$. The formula of a free wave in term of its radiation profile can also be given explicitly, see \cite{shenradiation}, for example. In this work we focus on the 3D radial case: 
\[
 u(r,t) = \frac{1}{r} \int_{t-r}^{t+r} G_-(s) {\rm d}s.  
\]
A basic calculation gives the initial data in term of the radiation profile 
\begin{align} \label{initial data by radiation profile} 
 &u_0(r) =  \frac{1}{r} \int_{-r}^{r} G_-(s) {\rm d}s; & & u_1(r) = \frac{G_-(r)-G_-(-r)}{r}.
\end{align}
The following relationship between the radiation profiles and the energy in the exterior region is useful in further argument. 
\begin{lemma} \label{tail G}
 Let $(u_0,u_1) \in \dot{H}^1\times L^2(\Rm^3)$ be radial initial data, whose radiation profile in the negative time direction is $G_-(s)$. Then we have
 \[
  \|(u_0,u_1)\|_{\mathcal{H}(R)}^2 = 8\pi \|G_-\|_{L^2(\{s: |s|>R\})}^2 + 4\pi R |u_0(R)|^2. 
 \]
\end{lemma} 
\begin{proof}
As a direct consequence of \eqref{initial data by radiation profile}, we have  
\[
 \int_{R}^\infty \left(\left|\partial_r (r u_0(r))\right|^2 + |ru_1(r)|^2 \right) {\rm d} r = 2 \int_R^\infty \left(|G_-(r)|^2 + |G_-(-r)|^2\right) {\rm d} r = 2\|G_-\|_{L^2(\{s:|s|>R\})}^2. 
\]
Next we apply integration by parts and obtain 
\begin{align*}
 \int_{R}^\infty \left|\partial_r (r u_0(r))\right|^2 {\rm d} r & = \int_R^\infty \left(r^2 |\partial_r u_0(r)|^2 + r \partial_r (|u_0(r)|^2) + |u_0(r)|^2 \right) {\rm d} r\\
 & = \int_R^\infty |\partial_r u_0(r)| r^2 {\rm d} r - R |u_0(R)|^2. 
\end{align*}
A combination of the identities above yields 
\[
 \int_{R}^\infty \left(\left|\partial_r u_0(r)\right|^2 + |u_1(r)|^2 \right) r^2 {\rm d} r = 2\|G_-\|_{L^2(\{s:|s|>R\})}^2 + R |u_0(R)|^2. 
\]
A change of variables then gives the desired result. 
\end{proof}

\begin{remark} \label{double tail G}
 A direct consequence of Lemma \ref{tail G} is 
 \[
  \|(u_0,u_1)\|_{\mathcal{H}(R)}^2 \geq 8\pi \|G_-\|_{L^2(\{s: |s|>R\})}^2, \qquad \forall R > 0. 
 \]
 It immediately follows that if $0\leq R_1 < R_2$, then 
 \[
  \|(u_0,u_1)\|_{\mathcal{H}(R_1)}^2 \leq \|(u_0,u_1)\|_{\mathcal{H}(R_2)}^2 + 8\pi \|G_-\|_{L^2(\{s: R_1 < |s|<R_2\})}^2 + 4\pi R_1 |u_0(R_1)|^2. 
 \]
\end{remark}

\subsection{Nonlinear radiation profiles} \label{sec: nonlinear profiles}

\begin{lemma} [Radiation fields of inhomogeneous equation] \label{scatter profile of nonlinear solution}
 Assume that $R\geq 0$. Let $u$ be a radial exterior solution to the wave equation
 \[
  \left\{\begin{array}{ll} \partial_t^2 u - \Delta u = F(t,x); & (x,t)\in \Omega_R; \\
  (u,u_t)|_{t=0} = (u_0,u_1) \in \dot{H}^1\times L^2. & \end{array} \right.
 \]
 If $F$ is a radial function satisfying $\|\chi_R F\|_{L^1 L^2(\Rm\times \Rm^3)}< +\infty$, then there exist unique radiation profiles $G_\pm \in L^2([R,+\infty))$ such that 
 \begin{align}
  \lim_{t\rightarrow +\infty} \int_{R+t}^\infty \left(\left|G_+(r-t) - r u_t (r, t)\right|^2 + \left|G_+(r-t) + r u_r (r, t)\right|^2\right) {\rm d}r & = 0; \label{positive ra}\\
  \lim_{t\rightarrow -\infty} \int_{R-t}^\infty \left(\left|G_-(r+t) - r u_t(r,t)\right|^2 +  \left|G_-(r+t) - r u_r(r,t)\right|^2\right) {\rm d} r & = 0. \label{negative ra}
 \end{align}
 In addition, the following estimates hold for $G_\pm$ given above and the corresponding radiation profiles $G_{0,\pm}$ of the initial data $(u_0,u_1)$:
 \begin{align*}
 2\sqrt{2\pi} \|G_- - G_{0,-}\|_{L^2([R,+\infty))} & \leq \|\chi_R F\|_{L^1 L^2((-\infty,0]\times \Rm^3)}; \\
  2\sqrt{2\pi}  \|G_+ - G_{0,+}\|_{L^2([R,+\infty))} & \leq \|\chi_R F\|_{L^1 L^2([0,+\infty)\times \Rm^3)}.
 \end{align*}
\end{lemma}
\begin{proof}
 The proof of a similar result has been given in the author's previous work \cite{ecarbitrary}. We still sketch the proof here for the reason of completeness. We may extend the domain of $u$ to the whole space-time $\Rm^3 \times \Rm$ by defining
 \[
  u = \mathbf{S}_L (u_0,u_1) + \int_0^t \frac{\sin (t-t')\sqrt{-\Delta}}{\sqrt{-\Delta}} [\chi_R(\cdot,t') F(\cdot,t')] {\rm d} t'. 
 \]
 In other words, the solution $u$ solves the wave equation $\partial_t^2 u - \Delta u = \chi_R F$ in the whole space-time. Since $\chi_R F \in L^1 L^2(\Rm\times \Rm^3)$, the solution $u$ must scatter in both two time directions, i.e. there exist two finite-energy free waves $u^\pm$, such that 
 \begin{equation} \label{scattering of u ex} 
  \lim_{t\rightarrow \pm \infty} \left\|\vec{u}(t) - \vec{u}^\pm (t)\right\|_{\dot{H}^1\times L^2} = 0. 
 \end{equation} 
 We let $G_+ \in L^2(\Rm)$ be the radiation profile of $u^+$ in the positive time direction and $G_-\in L^2(\Rm)$ be the radiation profile of $u^-$ in the negative time direction. Thus 
  \begin{align*}
  \lim_{t\rightarrow +\infty} \int_{0}^\infty \left(\left|G_+(r-t) - r u_t^+ (r, t)\right|^2 + \left|G_+(r-t) + r u_r^+ (r, t)\right|^2\right) {\rm d}r & = 0; \\
  \lim_{t\rightarrow -\infty} \int_{0}^\infty \left(\left|G_-(r+t) - r u_t^-(r,t)\right|^2 +  \left|G_-(r+t) - r u_r^-(r,t)\right|^2\right) {\rm d} r & = 0. 
 \end{align*}
 A combination of these two limits with \eqref{scattering of u ex} immediately yields \eqref{positive ra} and \eqref{negative ra}. If $G_+, \tilde{G}_+ \in L^2([R,+\infty))$ both satisfy \eqref{positive ra}, then  
 \[
  \int_R^\infty |G_+(s) - \tilde{G}_+(s)|^2 {\rm d} s = \lim_{t\rightarrow +\infty} \int_{R+t}^\infty |G_+(r-t) - \tilde{G}_+(r-t)|^2 {\rm d} r = 0. 
 \]
  This verifies the uniqueness of radiation profiles. Finally we verify the upper bound estimate of $\|G_\pm - G_{0,\pm}\|_{L^2([R,+\infty)}$. For convenience we introduce the following notations
 \begin{align*}
  &u^L = \mathbf{S}_L (u_0,u_1);& &v =  \int_0^t \frac{\sin (t-t')\sqrt{-\Delta}}{\sqrt{-\Delta}} [\chi_R(\cdot,t') F(\cdot,t')] {\rm d} t'.
 \end{align*}
 By $u = u^L + v$, we have
 \begin{align*}
  8\pi \int_{R}^\infty |G_+ (s) - G_{0,+}(s)|^2 {\rm d} s & = \lim_{t\rightarrow +\infty} 4\pi \int_{R+t}^\infty \left(\left|r u_t - r u_t^L\right|^2  + \left|r u_r - r u_r^L\right|^2 \right) {\rm d} r \\
  & = \lim_{t\rightarrow +\infty} \int_{|x|>R+t} |\nabla_{t,x} (u-u^L) |^2 {\rm d} x \\
  & \leq \lim_{t\rightarrow +\infty} \left\|\vec{u}(t) - \vec{u}^L(t)\right\|_{\dot{H}^1 \times L^2}^2 \\
  & = \lim_{t\rightarrow +\infty} \left\|\vec{v}(t)\right\|_{\dot{H}^1\times L^2}^2\\
  & \leq \|\chi_R F\|_{L^1 L^2(\Rm^+ \times \Rm^3)}^2. 
 \end{align*} 
 The negative time direction can be dealt with in the same manner. 
\end{proof}

\begin{remark} \label{remark channel}
 Let $u$ and $R$ be as above. The nonlinear radiation profile can be given by an explicit formula. In fact, the method of characteristic lines gives that 
 \[
  \frac{\rm d}{{\rm d} t} \left[(ru)_t (t+s,t) - (ru)_{r} (t+s,t)\right] = (s+t) F(s+t, t). 
 \]
 Without loss of generality we work as though the solutions are sufficiently smooth. By the theory of radiation fields, we have the $L^2$ convergence 
 \[ 
  (ru)_t (t+s,t) - (ru)_{r} (t+s,t) \rightarrow 2 G_+(s)
 \]
 This gives the formula 
 \[
  G_+ (s) - G_{0,+} (s) = \frac{1}{2} \int_{0}^\infty (s+t) F(s+t,t) {\rm d} t. 
 \]
 It immediately follows that for any $R' > R$, we have 
  \begin{align*}
 4\sqrt{\pi} \|G_- - G_{0,-}\|_{L^2(R,R')} & \leq \|\chi_{R,R'} F\|_{L^1 L^2((-\infty,0]\times \Rm^3)}; \\
  4\sqrt{\pi}  \|G_+ - G_{0,+}\|_{L^2(R,R')} & \leq \|\chi_{R,R'} F\|_{L^1 L^2([0,+\infty)\times \Rm^3)}.
 \end{align*}
 The constant here is better than Lemma \ref{scatter profile of nonlinear solution}, because a significant portion of energy for $\vec{v}$ is actually located in the cone $\{(x,t): |x|<|t|+R\}$, which we simply ignore in the proof of Lemma \ref{scatter profile of nonlinear solution}. 
\end{remark}

\subsection{Asymptotically equivalent solutions}

Assume that $u, v \in \mathcal{C}(\Rm; \dot{H}^1\times L^2)$. We say that $u$ and $v$ are $R$-weakly asymptotically equivalent if
\[
 \lim_{t\rightarrow \pm \infty} \int_{|x|>R+|t|} |\nabla_{t,x} (u-v)|^2 {\rm d} x = 0.
\]
Here $R\geq 0$ is a constant. In particular, if $R=0$, then we say $u$ and $v$ are asymptotically equivalent to each other. Since the integral above only involves the values of $u, v$ in the exterior region, the definition above also applies to exterior solutions. A solution $u$ is called ($R$-weakly) non-radiative solution if and only if it is asymptotically equivalent to zero. Non-radiative solutions, which play an essential role in the channel of energy method, have been extensively studied in recent years. Let us consider two examples. We start by considering an $R$-weakly non-radiative radial free wave $u$. It is equivalent to saying that the radiation profiles $G_\pm (s) = 0$ for $s>R$, or $G_-(s) = 0$ for $|s|>R$. An application of the explicit formula of linear free wave in term of radiation profile immediately gives 
\[
 u(r,t) = \frac{1}{r} \int_{-R}^R G_-(s) {\rm d} s, \qquad r>|t|+R. 
\]
This is a one-dimensional linear space spanned by $1/r$. Next we consider all radial non-radiative solutions to (CP1). One specific example of non-radiative solution is exactly the ground state mentioned in the introduction section
\[
 W(x) = \left(\frac{1}{3} + |x|^2\right)^{-1/2}. 
\] 
In fact this is the unique non-trivial non-radiative solution up to a rescaling/sign symmetry. In other words, all non-trivial non-radiative radial solutions can be given by 
\begin{align*}
 &\pm W_\lambda (x);& &W_\lambda (x) = \frac{1}{\lambda^{1/2}} W\left(\frac{x}{\lambda}\right), \qquad \lambda > 0. 
\end{align*}
We may also write them in the following form 
\[
 W^\alpha (x) = \frac{1}{\alpha} W\left(\frac{x}{\alpha^2}\right) = \frac{1}{\alpha} \left(\frac{1}{3} + \frac{|x|^2}{\alpha^4}\right)^{-1/2}, \qquad \alpha \in \Rm\setminus \{0\}. 
\]
which satisfies $W^\alpha (x) \simeq \alpha |x|^{-1}$ when $|x|$ is large. The notation $W^\alpha$ encodes both the scaling parameter $\lambda = \alpha^2$ and the sign into a single parameter $\alpha$. We will use both notations $W^\alpha$ and $W_\lambda$ frequently in the argument of this paper for convenience. 

In this work we need to consider asymptotically equivalent solutions of nonzero linear free waves. This generalize the conception of non-radiative solutions. In fact, thanks to the finite speed of propagation and a centre cut-off technique, we may show that any solution to (CP1) is $R$-weakly asymptotically equivalent to a linear free wave as long as $R$ is sufficiently large, by extending the domain of solution if necessary. If $u$ is a radial exterior solution to (CP1) defined in the exterior region $\Omega_R$, we may give a sufficient and necessary condition for $u$ to be $R$-weakly asymptotically equivalent to some linear free wave.
\begin{lemma} \label{finite eq scattering}
 Let $u$ be a radial exterior solution to (CP1) defined in $\Omega_R$, then $u$ is $R$-weakly asymptotically equivalent to some finite-energy linear free wave $w_L$ if and only if $\|\chi_R u\|_{Y(\Rm)} < +\infty$. 
\end{lemma}
\begin{proof}
 If $\|\chi_R u\|_{Y(\Rm)} < +\infty$, then we have $\|\chi_R F(u)\|_{L^1 L^2(\Rm \times \Rm^3)} < +\infty$. The existence of asymptotically equivalent free waves has been given in the proof of Lemma \ref{scatter profile of nonlinear solution}. Conversely if $u$ is $R$-weakly asymptotically equivalent to a free wave $w_L$, then we have 
 \[
  \lim_{t\rightarrow +\infty} \left\|\vec{w}_L(t) - \vec{u}(t)\right\|_{\mathcal{H}(R+t)} = 0. 
 \]
 A combination of this limit with the finite speed of propagation and the fact 
 \[
  \lim_{t \rightarrow +\infty} \|w_L\|_{Y([t,+\infty))} = 0
 \]
 yields 
 \[
  \lim_{t\rightarrow +\infty} \|\chi_{R+t} \mathbf{S}_L(\vec{u}(t))\|_{Y(\Rm^+)} = 0. 
 \]
 The small data theory and the uniqueness of exterior solution then guarantees that 
 \[
  \|\chi_R u \|_{Y([t,+\infty))} < +\infty, \qquad \forall t \gg 1. 
 \]
 This yields the estimate $\|\chi_R u\|_{Y(\Rm^+)}<+\infty$, since the definition of exterior solutions guarantees that $\|\chi_R u\|_{Y([0,t])}<+\infty$ for all $t>0$. The negative time direction can be dealt with in a similar way. 
\end{proof}

Given a radial finite-energy free wave $v_L$, all the possible radial weakly asymptotically equivalent solutions of $v_L$ have been discussed in Shen \cite{ecarbitrary}. Please note that the following result holds for all energy critical nonlinear terms $F(u)$ satisfying 
\begin{align*}
 &F(0) = 0;& & |F(u) - F(v)| \lesssim |u-v|(|u|^4 + |v|^4). 
\end{align*}

\begin{theorem}[One-parameter family, see \cite{ecarbitrary}]  \label{thm one parameter}
Let $w_L$ be a finite-energy radial free wave. Then there exists a one-parameter family $\{(u^\alpha, R_\alpha)\}_{\alpha \in \Rm}$ so that each pair $(u^\alpha, R_\alpha)$ satisfies either of the following
\begin{itemize}
 \item[(a)] The radial function $u^\alpha$ is an exterior solution to (CP1) in $\Omega_0$ and is asymptotically equivalent to $w_L$. In this case we choose $R_\alpha = 0^-$;
 \item[(b)] The radial function $u^\alpha$ is defined in $\Omega_{R_\alpha}$ with $R_\alpha \geq 0$ and $\|\chi_{R_\alpha} u^\alpha\|_{Y(\Rm)} = +\infty$ such  that for any $R > R_\alpha$, $u^\alpha$ is an exterior solution to (CP1) in $\Omega_R$ and is $R$-weakly asymptotically equivalent to $w_L$.
\end{itemize}
In addition, if $u$ is a radial exterior solution to (CP1) defined in $\Omega_R$ such that $u$ is $R$-weakly asymptotically equivalent to $w_L$, then there exists a unique real number $\alpha$, such that\footnote{In the case (a), we understand $0>R_\alpha = 0^-$.} $R>R_\alpha$ and $u(x,t) = u^\alpha(x,t)$ for $(x,t) \in \Omega_R$. We call the number $\alpha$ the characteristic number of $u$. The characteristic number can also be characterized by the asymptotic behaviour of the solutions. More precisely, given $\alpha, \beta \in \Rm$, we have 
\[
 \lim_{r\rightarrow +\infty} r^{1/2} \sup_{t\in \Rm} \left\|\vec{u^\alpha} (\cdot,t) - \vec{u^\beta} (\cdot, t) - ((\alpha-\beta)|x|^{-1}, 0)\right\|_{\mathcal{H}(|t|+r)} = 0. 
\]
\end{theorem}
\begin{remark}
 The main result of this work can be proved without using Theorem \ref{thm one parameter}. Please see Remark \ref{avoid one parameter}. We still introduce this conception of one parameter family for reason of completeness. 
\end{remark}

\section{Energy norm estimates of global solutions}

In this section we discuss the upper bound of the least energy norm in a long time interval. The main result of this section is 
\begin{lemma} \label{local boundedness}
There exists a small constant $\varepsilon_0 > 0$ and a large constant $K_0\gg 1$, such that if 
\begin{itemize}
 \item $u$ is a radial solution to (CP1) defined in a maximal time interval $(-T_-,T_+)$ with an energy $E>0$;
 \item The initial data $(u_0,u_1)$ satisfy $\|(u_0,u_1)\|_{\mathcal{H}(R)} < \varepsilon_0$;  
\end{itemize} 
then for any time $R \leq T< T_+/3$, there exists a time $t \in [T, 3T]$ satisfying
\[
 \|\vec{u}(\cdot,t)\|_{\dot{H}^1 \times L^2}^2 \leq 6E+K_0. 
\]
\end{lemma}
\begin{remark}
 The proof is based on the virial identity. This argument dates back to Levine \cite{negativeenergy}. Levine showed that any solution with a negative energy must blow up in finite time. It was proved in Duyskaerts-Kenig-Merle \cite{se} by a similar argument that if $T_+ = +\infty$, then 
 \[
  \liminf_{t\rightarrow+\infty} \|\vec{u}(t)\|_{\dot{H}^1\times L^2}^2 \leq 3 E. 
 \]
 The upper bound given in Lemma \ref{local boundedness} is larger but applies to a finite (but long) time interval. Please note that this argument does not depend on the radial assumption. 
\end{remark}
\begin{proof}
We assume that $\|\vec{u}(t)\|_{\dot{H}^1 \times L^2}^2 > 6 E+K_0$ for all $t\in [T,3T]$ and deduce a contradiction. Let $\varphi: \Rm \rightarrow [0,1]$ be a smooth cut-off function satisfying 
 \[
  \varphi(s) = \left\{\begin{array}{ll} 1, & s\leq 2\\ 0, & s\geq 3. \end{array}\right.
 \]
 and $\phi(s) = \varphi^2(s)$. We then define 
 \[
  J(t) = \int_{\Rm^3} |u(x,t)|^2 \phi(|x|/t) {\rm d} x.
 \]
 A straight-forward calculation yields 
 \[
  J'(t) = 2 \int_{\Rm^3} u u_t \phi(|x|/t) {\rm d} x - \int_{\Rm^3} |u|^2 \phi'(|x|/t) \frac{|x|}{t^2} {\rm d} x; 
 \]
 and
 \begin{align*}
  J''(t) & = 2 \int_{\Rm^3} |u_t|^2 \phi(|x|/t) {\rm d} x + 2 \int_{\Rm^3} u u_{tt} \phi(|x|/t) {\rm d} x - 4 \int_{\Rm^3} u u_t \phi'(|x|/t) \frac{|x|}{t^2} {\rm d} x \\
  & \qquad + \int_{\Rm^3} |u|^2 \phi''(|x|/t) \frac{|x|^2}{t^4} {\rm d} x + 2\int_{\Rm^3} |u|^2 \phi'(|x|/t) \frac{|x|}{t^3} {\rm d} x. 
 \end{align*}
 Inserting the equation $u_{tt} = \Delta u + |u|^4 u$ and integrating by parts, we obtain 
 \begin{align*}
  J''(t) & = 2 \int_{\Rm^3} (|u_t|^2-|\nabla u|^2 +|u|^6) \phi(|x|/t) {\rm d} x - 2\int_{\Rm^3} \phi'(|x|/t)u \nabla u \cdot \frac{x}{|x|t} {\rm d} x \\
  & \quad - 4 \int_{\Rm^3} u u_t \phi'(|x|/t) \frac{|x|}{t^2} {\rm d} x + \int_{\Rm^3} |u|^2 \phi''(|x|/t) \frac{|x|^2}{t^4} {\rm d} x + 2\int_{\Rm^3} |u|^2 \phi'(|x|/t) \frac{|x|}{t^3} {\rm d} x.
 \end{align*}
 By the finite speed of propagation, the small data theory, Hardy's inequality, we have 
 \begin{equation} \label{exterior estimate ft}
  \int_{|x|>|t|+R} \left(|\nabla u(x,t)|^2 + |u_t(x,t)|^2 + \frac{|u(x,t)|^2}{|x|^2} + |u(x,t)|^6\right) {\rm d} x \lesssim_1 \varepsilon_0^2. 
 \end{equation}
 Combining this with the facts 
 \begin{itemize}  
  \item $\phi(|x|/t)-1$ is nonzero only for $|x|>2 t \geq t+R$;
  \item $\phi'(|x|/t)$ and $\phi''(|x|/t)$ are nonzero only for $t+R\leq 2t < |x| < 3t$;
 \end{itemize} 
 we may write 
 \begin{align}
  J''(t) & =  \int_{\Rm^3} \left(2|u_t|^2 - 2|\nabla u|^2 + 2|u|^6 \right) {\rm d} x  + O(\varepsilon_0^2) \nonumber\\
  & = \int_{\Rm^3} \left(6 |u_t|^2 + 2|\nabla u|^2\right) {\rm d} x + 2\|\vec{u}(t)\|_{\dot{H}^1 \times L^2}^2 - 12 E + O(\varepsilon_0^2). \label{expression of Jpp}
 \end{align}
 Here the error term $O(\varepsilon_0^2)$ satisfies $|O(\varepsilon_0^2)| \lesssim_1 \varepsilon_0^2$. As a result, if $\varepsilon_0$ is sufficiently small, we have 
 \begin{align}
  |J'(t)|^2 & \leq 5 \left(\int_{\Rm^3} u u_t \phi(|x|/t) {\rm d} x\right)^2 + 5 \left(\int_{\Rm^3} |u|^2 \phi'(|x|/t) \frac{|x|}{t^2} {\rm d} x\right)^2 \nonumber\\
  & \leq 5 \left(\int_{\Rm^3} |u_t|^2 {\rm d} x\right)\left(\int_{\Rm^3} |u|^2 \phi^2 (|x|/t) {\rm d} x\right) + \left(\int_{\Rm^3} |u|^2 \phi(|x|/t) {\rm d} x\right) O(\varepsilon_0^2) \nonumber\\
  & \leq \frac{5}{6}J''(t) J(t). \label{JJpJpp}
 \end{align}
 It is not difficult to see that $J(t) \in \mathcal{C}^2([T,3T])$. Now we assume that $J''(t)$ takes its minimum value $M_0$ at time $t_0$ in the time interval $[T,3T]$. By the expression \eqref{expression of Jpp} of $J''(t)$ and our assumption $\|\vec{u}(t)\|_{\dot{H}^1\times L^2}^2 > 6 E + K_0$, we have 
 \[
  M_0 \geq \int_{\Rm^3} \left(6 |u_t(x,t_0)|^2 + 2|\nabla u(x,t_0)|^2\right) {\rm d} x \geq 2\|\vec{u}(t_0)\|_{\dot{H}^1\times L^2}^2 \geq 12 E + 2 K_0. 
 \]
 By the energy conservation law and the assumption $E>0$, we also have 
 \[
  \int_{\Rm^3} |u(x,t_0)|^6 {\rm d} x \leq 3 \|\vec{u}(t_0)\|_{\dot{H}^1 \times L^2}^2 \leq 3M_0/2. 
 \]
 By H\"{o}lder inequality and \eqref{JJpJpp}, we have 
 \begin{align*}
  J(t_0) & \leq \int_{|x|<3t_0} |u(x,t_0)|^2 {\rm d} x \lesssim_1 T^2 M_0^{1/3}; \quad \Longrightarrow \quad |J'(t_0)| \lesssim_1 T M_0^{2/3}. 
 \end{align*}
 There are two cases: Case one, $t_0 \leq 2T$. We have 
 \begin{align*}
   J'(t_0 + T/10) & = J'(t_0) + \int_{t_0}^{t_0+T/10} J''(t) {\rm d} t;\\
   J(t_0 + T/10) & = J(t_0) + \frac{T}{10} J'(t_0) + \int_{t_0}^{t_0+T/10} \left(t_0+T/10-t\right) J''(t) {\rm d} t.
 \end{align*}
 Since $J''(t) \geq M_0 > 0$ for $t\in [t_0,t_0+T/10]\subset [T,3T]$, the following inequalities hold:
 \begin{align*}
  \int_{t_0}^{t_0+T/10} J''(t) {\rm d} t & \gtrsim_1 M_0 T; \\
  \int_{t_0}^{t_0+T/10} \left(t_0+T/10-t\right) J''(t) {\rm d} t & \gtrsim_1 T^2 M_0. 
 \end{align*}
 This implies that if $K_0$, thus $M_0$ is sufficiently large, then the integral part is the dominating term in the expressions of $J(t_0+T/10)$ and $J'(t_0+T/10)$ given above. In addition, it is clear that 
 \[
  \int_{t_0}^{t_0+T/10} J''(t) {\rm d} t  \geq \frac{10}{T}\int_{t_0}^{t_0+T/10} \left(t_0+T/10-t\right) J''(t) {\rm d} t.
 \]
 Therefore the following inequalities hold as long as the constant $K_0$ is sufficiently large:
 \begin{align*}
  &J(t_0+T/1)), J'(t_0+T/10) > 0;& & \frac{J'(t_0 + T/10)}{J(t_0 + T/10)} \geq \frac{9}{T}.
 \end{align*}
 We define $Q(t) = J'(t) /J(t) > 0$ for all $t\in [t_0+T/10,3T]$. The inequality \eqref{JJpJpp} then gives
 \[
  Q'(t) = \frac{J''(t) J(t) - (J'(t))^2}{J(t)^2} \geq \frac{1}{5} \left(\frac{J'(t)^2}{J(t)^2}\right) = \frac{1}{5} Q^2(t), \qquad t\in  [t_0+T/10,3T];
 \]
 which implies that 
 \[
  \frac{1}{Q(t_0+T/10)} - \frac{1}{Q(3T)} \geq \frac{1}{5} \left[3T - (t_0+T/10)\right] \geq \frac{9}{50} T.
 \]
 This is a contradiction since $Q(t_0+T/10) \geq 9/T$ and $Q(3T)>0$. Now let us consider case two, namely $t_0 \in (2T,3T]$. In this case we consider 
 \begin{align*}
   J'(t_0 - T/10) & = J'(t_0) - \int_{t_0-T/10}^{t_0} J''(t) {\rm d} t;\\
   J(t_0 - T/10) & = J(t_0) - \frac{T}{10} J'(t_0) + \int_{t_0-T/10}^{t_0} \left(t-t_0+T/10\right) J''(t) {\rm d} t.
 \end{align*}
 A similar argument gives
 \begin{align*}
  &J(t_0-T/10) > 0, \; J'(t_0-T/10) < 0,& &|Q(t_0 - T/10)| = \frac{|J'(t_0 - T/10)|}{|J(t_0 - T/10)|} \geq \frac{9}{T}, 
 \end{align*}
 which gives a contradiction if we consider the evolution of $Q(t)$ as the time decreases from $t_0-T/10$ to $T$. 
 \end{proof}

\section{Soliton resolution of almost non-radiative solutions} 
The following proposition separate each bubble one-by-one as long as the radiation is sufficiently weak in the main light cone. This is the most important observation in this work. 
\begin{proposition} \label{main tool} 
 Let $n$ be a positive integer. Then there exists a small constant $\delta_0 = \delta_0(n)>0$ and an absolute constant $c_2 \gg 1$, such that if $v_L$ is a finite-energy radial free wave with $\delta \doteq \|\chi_0 v_L\|_{Y(\Rm)} < \delta_0$, then any weakly asymptotically equivalent solution $u$ to (CP1) of $v_L$ satisfies either of the following: (we extend the domain of $u$ if necessary)
 \begin{itemize}
  \item [(a)] The solution $u$ is an exterior solution in $\Omega_0$. In addition, there exists a sequence $\{\alpha_j\}_{j=1,2,\cdots,J}$ with $0 \leq J\leq n-1$ and $|\alpha_1| > |\alpha_2| > \cdots > |\alpha_J|>0$ such that  
  \begin{align*}
    \frac{|\alpha_{j+1}|}{|\alpha_j|} & \lesssim_j \delta, \qquad j=1,2,\cdots, J-1; \\
   \left\|\vec{u}(\cdot, 0)-\sum_{j=1}^J (W^{\alpha_j},0) - \vec{v}_L(\cdot,0)\right\|_{\dot{H}^1\times L^2} & + \left\|\chi_0 \left(u - \sum_{j=1}^J W^{\alpha_j} \right)\right\|_{Y(\Rm)}  \lesssim_J \delta.
  \end{align*}
  \item[(b)] There exists a sequence $\{\alpha_j\}_{j=1,2,\cdots,n}$ with $|\alpha_1| > |\alpha_2| > \cdots > |\alpha_n|>0$ and 
  \[
     \frac{|\alpha_{j+1}|}{|\alpha_j|} \lesssim_j \delta, \qquad j=1,2,\cdots, J-1;
  \]
  such that $u$ is an exterior solution in the region $\Omega_{c_2 \alpha_n^2}$ and satisfies 
  \begin{align*}
    \left\|\vec{u}(\cdot,0)-\sum_{j=1}^n (W^{\alpha_j},0) - \vec{v}_L(\cdot,0)\right\|_{\mathcal{H}(c_2 \alpha_n^2)} + \left\|\chi_{c_2 \alpha_n^2} \left(u - \sum_{j=1}^n W^{\alpha_j} \right)\right\|_{Y(\Rm)} & \lesssim_n \delta.
  \end{align*}
 \end{itemize}
 \end{proposition} 

\begin{remark} \label{avoid one parameter}
 The domain extension we make above still guarantees that $u$ is (weakly) asymptotically equivalent to $v_L$ in the corresponding exterior region. If we recall the conception of one-parameter family given in Theorem \ref{thm one parameter}, Proposition \ref{main tool} claims that any solution $u^\alpha$ in the one-parameter family satisfies either (a) or (b). In addition, we must have $R_\alpha = 0^-$ in case (a); or $R_\alpha< c_2 \alpha_n^2$ in case (b). This is exactly the version proved in this section. Please note that $\|\chi_R u^\alpha\|_{Y(\Rm)} < +\infty$ holds for any $R > R_\alpha$ by Lemma \ref{finite eq scattering}. If we assume that $u$ is an exterior solution in $\Omega_0$ and asymptotically equivalent to $v_L$, then the same proof shows that the same conclusion of Proposition \ref{main tool} still holds without using the conception of one-parameter family or an extension of domain. It is exactly the case in the proof of our main result, i.e. Theorem \ref{main thm}. 
 \end{remark}

\begin{remark} \label{nonlinear estimate} 
 A direct calculation of nonlinear estimate shows that if $\delta < \delta_0(n)$ is sufficiently small, then a solution $u$ in case (a) satisfies 
\[
 \left\|\chi_{0} \left(F(u) - \sum_{j=1}^J F(W^{\alpha_j}) \right)\right\|_{L^1 L^2(\Rm \times \Rm^3)} \lesssim_{J} \delta.
\]
 Similarly a solution in case (b) satisfies 
\[
 \left\|\chi_{c_2 \alpha_n^2} \left(F(u) - \sum_{j=1}^n F(W^{\alpha_j}) \right)\right\|_{L^1 L^2(\Rm \times \Rm^3)} \lesssim_{n} \delta . 
\]
Here we use the estimates (the estimates for terms with three or more different scales follow from an interpolation)
\[
 \left\|\chi_0 (W^{\alpha})^k (W^{\alpha'})^{5-k}\right\|_{L^1 L^2(\Rm\times \Rm^3)} \lesssim_1 \min\left\{\frac{|\alpha|}{|\alpha'|}, \frac{|\alpha'|}{|\alpha|}\right\}, \qquad k=1,2,3,4. 
\]
Please note that in this work $\delta_0(n)$, $\delta_{0}(n,c)$ or similar notations represent small positive constants depending on $n$ (or $n$ and $c$). They may represent different constants at different places. 
\end{remark} 

\begin{remark}
 The estimate $|\alpha_{j+1}|/|\alpha_j| \lesssim_j \delta$ given in Proposition \ref{main tool} are optimal if we ignore a constant multiple. In fact, we may consider the exterior solution with initial data 
 \[
  \vec{u}(0) = (W^{\alpha_1} + W^{\alpha_2}, 0). 
 \]
 Here $|\alpha_2|/|\alpha_1| \ll 1$. Comparing this with the approximated solution $\tilde{u} = W^{\alpha_1} + W^{\alpha_2}$ and applying the perturbation theory, we obtain that the asymptotically equivalent free wave $v_L$ of $u$ satisfies 
 \[
  \|v_L\|_{Y(\Rm)} \lesssim_1 \|\vec{v}_L\|_{\mathcal{H}} \lesssim_1 \left\|\chi_0(F(W^{\alpha_1} + W^{\alpha_2})-F(W^{\alpha_1})-F(W^{\alpha_2}))\right\|_{L^1 L^2}\lesssim_1 |\alpha_2| /|\alpha_1|. 
 \]
 An application of Proposition \ref{main tool} gives 
 \begin{align*}
  &u(\cdot, 0) = W^{\tilde{\alpha}_1} + W^{\tilde{\alpha}_2} + v_L(\cdot,0) + O(\|\chi_0 v_L\|_{Y(\Rm)});& &|\tilde{\alpha}_2|/|\tilde{\alpha}_1| \lesssim_1 \|\chi_0 v_L\|_{Y(\Rm)}. 
 \end{align*}
 Although $\tilde{\alpha}_j$ is not necessary the same as $\alpha_j$, but we must have $\tilde{\alpha}_j/\alpha_j \simeq 1$. It follows that 
 \[
  |\alpha_2|/|\alpha_1| \simeq |\tilde{\alpha}_2|/|\tilde{\alpha}_1| \simeq \|\chi_0 v_L\|_{Y(\Rm)} \simeq \|\vec{v}_L\|_{\mathcal{H}}.  
 \]
\end{remark}

The rest of this section is devoted to the proof of Proposition \ref{main tool}. We will apply an induction in the positive integer $n$. More precisely we split the proof into three steps, each presented in an individual subsection (Section 4.2-4.4).
\begin{itemize}
 \item In Step one we show that Proposition \ref{main tool} holds for $n=1$; 
 \item In Step two we give some key estimates, which lead to the ratio estimate $|\alpha_{j+1}|/|\alpha_j| \lesssim_j \delta$; 
 \item In Step three we prove that the proposition holds for $n+1$ as long as the proposition holds for $n$. 
\end{itemize}
Before we start the proof, we first introduce some notations and preliminary results. 

\subsection{Preliminary results} \label{sec: pre}

\paragraph{Notations} We first introduce a few notations. Given a sequence $\{\alpha_j\}_{j=1,2,\cdots,n}$, we define 
\begin{align*}
 S_n(x,t) & = v_L + \sum_{j=1}^n W^{\alpha_j}; & e_n(x,t) = \sum_{j=1}^n F(W^{\alpha_j}) - F(S_n).
\end{align*}
$S_n$ solves the following wave equation in the region $\Omega_{0}$. 
\[
 (\partial_t^2 - \Delta) S_n = F(S_n) + e_n(x,t).
\]
We also define $w_n = u - S_n$ in an exterior region $\Omega_R$, as long as $u$ is a well-defined exterior solution in this region. Clearly $w_n$ is an exterior solution to 
\[
 (\partial_t^2 - \Delta) w_n = F(u) - F(S_n) - e_n(x,t).
\]
Let $w_{n,L}$ and $G_n$ be the linear free wave and radiation profiles with the initial data $\vec{w}_n(0)$. In this section all radiation profiles are the one in the negative time direction, unless specified otherwise. When there is no risk of confusion, we use notations $w, w_L, G$ respectively.  

\begin{remark}
 If $w$ is only defined in an exterior region $\Omega_R$ with $R>0$, then its initial data $\vec{w}(0)$ are not uniquely determined by $w$. However, $\vec{w}(0)$ are uniquely determined by $w$ in the exterior region $\{x: |x|>R\}$. This implies that $w_L$ are uniquely determined in the exterior region $\Omega_R$. In addition, the radiation profile $G$ are also uniquely determined in the space $L^2(\{s: |s|>R\})$. Although $G(s)$ can not be uniquely determined for $s\in (-R,R)$, the integral 
 \[
  \int_{-R}^R G(s) {\rm d} s
 \]
 is uniquely determined by $w$. These properties about $G$ immediately follows from formula \eqref{initial data by radiation profile}. For convenience we let $\Omega_{R'}$ be the maximal region of existence for the solution $u$ in this section. 
\end{remark}

\begin{lemma} \label{lemma inequality}
 Let $u$, $S$ be exterior solutions of (CP1) and $(\partial_t^2 - \Delta) S = F(S) + e(x,t)$ in $\Omega_{R_1}$, respectively, with 
 \begin{align*}
  \|\chi_{R_1} u\|_{Y(\Rm)}, \|\chi_{R_1} S\|_{Y(\Rm)}, \|\chi_{R_1} e(x,t)\|_{L^1 L^2(\Rm\times \Rm^3)} < +\infty.  
 \end{align*}
 Let $w = u - S$ and $G$ be the radiation profile of the initial data $\vec{w}(0)$. There exists an absolute constant $C_1\geq 1$ such that the following inequality holds for any $R_2 > R_1\geq 0$:
 \begin{align*}
  \|\chi_{R_1} w\|_{Y(\Rm)} &\leq C_1 \left(R_1^{1/2}|w(R_1,0)| + \|G\|_{L^2(\{s: R_1<|s|<R_2\})} + \|\vec{w}(\cdot,0)\|_{\mathcal{H}(R_2)}\right)\\
  & \qquad + C_1\left(\left\|\chi_{R_2}\left(F(u)-F(S) - e(x,t)\right)\right\|_{L^1 L^2(\Rm \times \Rm^3)} + \|\chi_{R_1,R_2}e(x,t)\|_{L^1 L^2(\Rm \times \Rm^3)}\right)\\
  & \qquad \qquad + C_1\left(\|\chi_{R_1,R_2} w\|_{Y(\Rm)}^4 + \|\chi_{R_1,R_2}S\|_{Y(\Rm)}^4\right)\|\chi_{R_1,R_2} w\|_{Y(\Rm)}.
 \end{align*}
  Here $\chi_{R_1,R_2}$ is the characteristic function of the region (as defined at the beginning of Section 2)
 \[
  \Omega_{R_1,R_2} = \{(x,t): |t|+R_1<|x|<|t|+R_2\}. 
 \]
 In addition, the inequality $|F(x+y)-F(y)|\leq C_1 |x|(|x|^4+|y|^4)$ holds for all numbers $x, y$. 
\end{lemma}
\begin{proof}
 It is sufficient to prove the first inequality, because the second inequality clearly holds for a sufficiently large constant $C_1$. First of all, we may apply Strichartz estimates, as well as Remark \ref{double tail G}, and obtain 
 \begin{align*}
  \|\chi_{R_1} w_L\|_{Y(\Rm)} & \lesssim_1 \|(w(0), w_t(0))\|_{\mathcal{H}(R_1)} \\
  & \lesssim_1 R_1^{1/2}|w(R_1,0)| + \|G\|_{L^2(\{s: R_1<|s|<R_2\})} + \|(w(0), w_t(0))\|_{\mathcal{H}(R_2)}. 
 \end{align*}
 Here $w_L$ is the linear free wave with initial data $\vec{w}(0)$. Since $w$ satisfies the equation $(\partial_t^2 - \Delta) w = F(u) - F(S) - e(x,t)$, we have 
 \begin{align*}
  \|\chi_{R_1} w\|_{Y(\Rm)} & \lesssim_1 \|\chi_{R_1} w_L\|_{Y(\Rm)} + \left\|\chi_{R_1}\left(F(u) - F(S) - e(x,t)\right)\right\|_{L^1 L^2} \\
  & \lesssim_1 \|\chi_{R_1} w_L\|_{Y(\Rm)} + \left\|\chi_{R_2}\left(F(u) - F(S) - e(x,t)\right)\right\|_{L^1 L^2} \\
  & \qquad + \left\|\chi_{R_1,R_2}\left(F(u) - F(S)\right)\right\|_{L^1 L^2} + \left\|\chi_{R_1,R_2}e(x,t)\right\|_{L^1 L^2}. 
 \end{align*}
 Finally H\"{o}lder inequality gives 
 \begin{align*}
  \left\|\chi_{R_1,R_2}\left(F(u) - F(S)\right)\right\|_{L^1 L^2} \lesssim_1 \left(\|\chi_{R_1,R_2} w\|_{Y(\Rm)}^4 + \|\chi_{R_1,R_2}S\|_{Y(\Rm)}^4\right)\|\chi_{R_1,R_2} w\|_{Y(\Rm)}. 
 \end{align*}
 A combination of these inequalities finishes the proof. 
\end{proof} 

\begin{lemma} \label{lemma connection} 
 There exist absolute positive constants $\varepsilon_1$, $\beta$, $\eta$ such that if $0\leq R_1<R_2$ and
 \begin{itemize} 
  \item $u$ is an exterior solution to (CP1) and $S$ is an exterior solution to the equation 
 \[
  (\partial_t^2 -\Delta) S= F(S) + e(x,t),
 \] 
 both in the region $\Omega_{R_1}$, with $\|\chi_{R_1} u\|_{Y(\Rm)}, \|\chi_{R_1} S\|_{Y(\Rm)}, \|\chi_{R_1} e(x,t)\|_{L^1 L^2(\Rm\times \Rm^3)} < +\infty$; 
 \item both $u$, $S$ are asymptotically equivalent to each other in $\Omega_{R_1}$; 
 \item $u$, $S$ and $w=u-S$ satisfy the following inequalities
 \begin{align*}
  \varepsilon \doteq \|\vec{w}(\cdot,0)\|_{\mathcal{H}(R_2)}  + \|\chi_{R_1,R_2} e(x,t)\|_{L^1 L^2(\Rm \times \Rm^3)} \qquad & \\
   + \|\chi_{R_2}\left(F(u) - F(S)- e(x,t)\right)\|_{L^1 L^2(\Rm \times \Rm^3)} & \leq \varepsilon_1; \\
   \|\chi_{R_1,R_2} S\|_{Y(\Rm)} & \leq \eta; \\
  \sup_{R_1\leq r\leq R_2} \left(r^{1/2} |w(r,0)|\right) & \leq \beta;
 \end{align*}
 \end{itemize}
 then we have 
 \begin{align*}
   \left\|\vec{w}(\cdot,0)\right\|_{\mathcal{H}(R_1)} + \|\chi_{R_1} w\|_{Y(\Rm)} & \lesssim_1 R_1^{1/2}|w(R_1,0)| +\varepsilon.
 \end{align*}
\end{lemma}
\begin{proof}
 Let $w_L$ and $G$ be the linear free wave and radiation profile with initial data $\vec{w}(\cdot,0)$. By Lemma \ref{lemma inequality}, we obtain for any $R\in [R_1,R_2)$ that
 \begin{align*}
  \|\chi_{R} w\|_{Y(\Rm)} \leq C_1 \left(R^{1/2}|w(R,0)| + \|G\|_{L^2(\{s: R<|s|<R_2\})} + \|\chi_{R,R_2} w\|_{Y(\Rm)}^5+ \eta^4 \|\chi_{R,R_2} w\|_{Y(\Rm)} +\varepsilon \right).
 \end{align*}
 We choose $\eta$ to be a sufficiently small number such that $C_1 \eta^4 < 1/(4C_1) < 1/2$, thus 
  \begin{align} \label{recurrence w Y} 
  \|\chi_{R,R_2} w\|_{Y(\Rm)} \leq 2C_1 \left(R^{1/2}|w(R,0)| + \|G\|_{L^2(\{s: R<|s|<R_2\})} + \|\chi_{R,R_2} w\|_{Y(\Rm)}^5+\varepsilon \right).
 \end{align}
 We choose small constants $\varepsilon_1 = \beta$ such that 
 \begin{align*}
  2C_1 (8C_1\beta)^4 < \frac{1}{4C_1} < \frac{1}{4} \quad \Longrightarrow \quad  8C_1 \beta > 2C_1(3\beta + (8C_1\beta)^5). 
 \end{align*}
 As a result, if $\|G\|_{L^2(\{s: R<|s|<R_2\})} \leq \beta$, then a continuity argument in $R$ shows that 
 \[
  \|\chi_{R,R_2} w\|_{Y(\Rm)} < 8 C_1 \beta. 
 \]
 Inserting this to \eqref{recurrence w Y} and using the choice of $\beta$, we obtain 
 \[
  \|\chi_{R,R_2} w\|_{Y(\Rm)} \leq 2C_1 \left(R^{1/2}|w(R,0)| + \|G\|_{L^2(\{s: R<|s|<R_2\})} +\varepsilon \right) + \frac{1}{4} \|\chi_{R,R_2} w\|_{Y(\Rm)},
 \]
 which implies 
 \[
   \|\chi_{R,R_2} w\|_{Y(\Rm)} \leq \frac{8}{3}C_1 \left(R^{1/2}|w(R,0)| + \|G\|_{L^2(\{s: R<|s|<R_2\})} +\varepsilon \right).
 \]
 An application of the nonlinear radiation profile shows (we apply Lemma \ref{scatter profile of nonlinear solution} on $w$ and recall the choice of $\beta, \eta$)
 \begin{align*}
 2\sqrt{2\pi} \|G\|_{L^2(\{s:|s|>R\})} & \leq \|\chi_{R} (F(u)-e(x,t)-F(S))\|_{L^1 L^2 (\Rm \times \Rm^3)}\nonumber \\
 & \leq \|\chi_{R,R_2} (F(w+S)-F(S))\|_{L^1 L^2(\Rm \times \Rm^3)} + \varepsilon \nonumber\\
 & \leq C_1\left(\|\chi_{R,R_2} w\|_{Y(\Rm)}^5 + \|\chi_{R,R_2} S\|_{Y(\Rm)}^4 \|\chi_{R,R_2} w\|_{Y(\Rm)}\right) + \varepsilon \\
 & \leq C_1 \left((8C_1 \beta)^4 + \eta^4\right) \|\chi_{R,R_2} w\|_{Y(\Rm)} + \varepsilon \nonumber \\
 & \leq \frac{3}{8C_1}\cdot \frac{8}{3}C_1 \left(R^{1/2}|w(R,0)| + \|G\|_{L^2(\{s: R<|s|<R_2\})} +\varepsilon \right) + \varepsilon\\
 & \leq R^{1/2}|w(R,0)| + \|G\|_{L^2(\{s: R<|s|<R_2\})} +2\varepsilon.
\end{align*}
This immediately gives 
\[
 \|G\|_{L^2(\{s:R<|s|<R_2\})} \leq \frac{1}{4}R^{1/2} |w(R,0)| + \frac{1}{2} \varepsilon \leq \frac{3}{4}\beta.
\]
A continuity argument in $R$ shows that $\|G\|_{L^2(\{s: R_1<|s|<R_2\})} \leq 3\beta/4$. Thus the inequalities above hold for all $R\in [R_1,R_2)$. A combination of these inequalities with Remark \ref{double tail G} finishes the proof. 
\end{proof}

\subsection{Step one} 

In this subsection we prove Proposition \ref{main tool} for $n=1$. We let $S=S_0 = v_L$ and $w=u-S$. Please note that in this case $e(x,t) = -F(v_L)$. Thus 
\[
 \|\chi_0 e(x,t)\|_{L^1 L^2(\Rm \times \Rm^3)} \leq \delta^5. 
\]
Let $\varepsilon_1$, $\beta$, $\eta$ be constants in Lemma \ref{lemma connection}. This is clear that 
\begin{align*}
 \varepsilon(R_2) \doteq \|\vec{w}(0)\|_{\mathcal{H}_{R_2}} + \|\chi_{R_2} F(u)\|_{L^1 L^2}
\end{align*}
satisfies the limit 
\[
 \lim_{R_2 \rightarrow +\infty} \varepsilon(R_2) = 0.
\]
Next we choose an absolute constant $c_2 \gg 1$ such that 
 \begin{itemize} 
  \item The inequality $\displaystyle \|\chi_{c_2} W\|_{Y(\Rm)} < \eta/3$ holds. 
  \item The inequality $\beta_1 \doteq c_2^{1/2} \left(\frac{1}{3} + c_2^2\right)^{-1/2} < \beta/2$ holds. 
 \end{itemize} 
Now we consider the function $r^{1/2} |w(r,0)|$ defined for all nonnegative radius $r > R'$, where $R'$ is determined by the maximal domain $\Omega_{R'}$ of u. By the point-wise decay of radial $\dot{H}^1$ functions, we have
\[
 \lim_{r\rightarrow +\infty} r^{1/2} |w(r,0)| = 0.
\]
Let us assume $\delta < \min\{\eta/6, \varepsilon_1^{1/5}\}$. There are two cases: 
\paragraph{Case One} We have
\[
 \sup_{r>R'} r^{1/2} |w(r,0)| < \beta. 
\]
Now we are able to apply Lemma \ref{lemma connection} for any $R_1$ slightly larger than $R'$ and any sufficiently large $R_2$ to conclude 
 \begin{align*}
  \left\|\vec{w}(\cdot,0)\right\|_{\mathcal{H}(R_1)} + \|\chi_{R_1} w\|_{Y(\Rm)}  \lesssim_1 R_1^{1/2}|w(R_1,0)| + \delta^5 + \varepsilon(R_2). 
 \end{align*}
 The norms $\|\chi_{R_1} w\|_{Y(\Rm)}$, thus the norms $\|\chi_{R_1} u\|_{Y(\Rm)}$ are uniformly bounded for all $R_1>R'$. As a result of Theorem \ref{thm one parameter}, we must have $R' = 0^-$. This means that we may choose $R_1 = 0$ and obtain 
 \[
  \left\|\vec{w}(\cdot,0)\right\|_{\dot{H}^1 \times L^2} + \|\chi_{0} w\|_{Y(\Rm)} \lesssim_1 \delta^5 + \varepsilon(R_2). 
 \]
 Letting $R_2 \rightarrow +\infty$ yields that 
 \[
  \|\vec{u}(0) - \vec{v}_L(0)\|_{\dot{H}^1\times L^2} +  \|\chi_{0} (u-v_L)\|_{Y(\Rm)} \lesssim_1 \delta^5.
 \] 
 This verifies case (a) with $J=0$. 
 \paragraph{Case Two} We have
 \[
 \sup_{r>R'} r^{1/2} |w(r,0)| \geq \beta. 
\]
 Combining this with the continuity and the limit at the infinity, we may always find a radius $R_1 > R'$ such that 
 \[
  \sup_{r>R_1} r^{1/2} |w(r,0)| = R_1^{1/2} |w(R_1,0)| = \beta_1. 
 \]
 Now we choose $\alpha_1 = \pm (c_2^{-1} R_1)^{1/2}$, where the sign is equal to that of $w(R_1,0)$. A basic calculation shows that 
 \[
  R_1^{1/2} w(R_1,0) = R_1^{1/2} W^{\alpha_1}(R_1) = \pm \beta_1. 
 \]
 Now we let $S_1 = v_L + W^{\alpha_1}$, $w_1 = u - S_1 = w - W^{\alpha_1}$ and $e_1 = F(W^{\alpha_1})  - F(S_1)$. They satisfy
 \begin{align*}
  \lim_{R_2\rightarrow +\infty} \left(\|\vec{w}_1(0)\|_{\mathcal{H}(R_2)} + \|\chi_{R_2} (F(u) - F(S_1) -e_1(x,t))\|_{L^1 L^2(\Rm \times \Rm^3)}\right) & = 0;\\
  \|\chi_0 e_1(x,t)\|_{L^1 L^2(\Rm \times \Rm^3)} & \lesssim_1 \delta;\\
  \sup_{r\geq R_1} r^{1/2} |w_1(r,0)| & \leq \beta; \\
  \|\chi_{R_1} S_1\| & \leq \delta + \eta/3;\\
  |w_1(R_1,0)| & = 0.
 \end{align*}
 As a result, if $\delta$ also satisfies $\delta < \delta_0(1)$, where $\delta_0(1)$ is a very small absolute constant, we may apply Lemma \ref{lemma connection} for large radius $R_2$ to conclude that 
 \begin{align*}
  \|\chi_{R_1} w_1\|_{Y(\Rm)} + \left\|\vec{w}_1 (0)\right\|_{\mathcal{H}(R_1)} & \lesssim_1  \delta + \|\vec{w}_1(0)\|_{\mathcal{H}(R_2)} + \|\chi_{R_2} (F(u) - F(S_1) -e_1(x,t))\|_{L^1 L^2}.
 \end{align*}
 Making $R_2 \rightarrow +\infty$ verifies that $u$ satisfies (b) thus finishes the proof. 

\subsection{Step two} 

In this subsection we assume that Proposition \ref{main tool} holds for a positive integer $n$ and prove a few estimates, which leads to the inequality $|\alpha_{n+1}|/|\alpha_n| \lesssim_n \delta$. We first give a lemma, which is a modified version of Lemma \ref{lemma connection}. 

\begin{lemma} \label{lemma connection 2} 
 Let $\eta$ be the constant in Lemma \ref{lemma connection}. There exists an absolute positive constant $\varepsilon_2$ such that if $3 R_2/4 \leq R_1<R_2$ and
 \begin{itemize} 
  \item $u$ is an exterior solution to (CP1) and $S$ is an exterior solution to the equation 
 \[
  (\partial_t^2 -\Delta) S= F(S) + e(x,t). 
 \] 
 with $\|\chi_{R_1} u\|_{Y(\Rm)}, \|\chi_{R_1} S\|_{Y(\Rm)}, \|\chi_{R_1} e(x,t)\|_{L^1 L^2} < +\infty$. 
 \item Solutions $u$, $S$ are asymptotically equivalent to each other in $\Omega_{R_1}$. 
 \item $u$, $S$ and $w=u-S$ satisfy the following inequalities
 \begin{align*}
  \varepsilon \doteq \|\vec{w}(\cdot,0)\|_{\mathcal{H}(R_2)}  + \|\chi_{R_1,R_2} e(x,t)\|_{L^1 L^2(\Rm \times \Rm^3)} \qquad  & \\
 + \|\chi_{R_2}\left(F(u) - F(S)- e(x,t)\right)\|_{L^1 L^2(\Rm \times \Rm^3)} & \leq \varepsilon_2; \\
   \|\chi_{R_1,R_2} S\|_{Y(\Rm)} & \leq \eta; 
 \end{align*}
 \end{itemize}
 Then we have 
 \begin{align*}
  \|\chi_{R_1} w\|_{Y(\Rm)} + \left\|(w(\cdot,0), w_t(\cdot,0))\right\|_{\mathcal{H}(R_1)} \lesssim_1 \varepsilon. 
 \end{align*}
\end{lemma}

\begin{proof}
 The proof is similar to Lemma \ref{lemma connection}. Let $w_L$ and $G$ be the linear free wave and radiation profile with initial data $\vec{w}(\cdot,0)$. First of all, for $R\in [R_1,R_2)$ we have that
 \begin{align*}
  R^{1/2} |w(R,0)| & = \left|R^{-1/2} \int_{-R}^{R} G(s) {\rm d} s\right| \\
  & \leq R^{-1/2} \left|\int_{-R_2}^{R_2} G(s) {\rm d} s\right| + R^{-1/2} \int_{R<|s|<R_2} |G(s)| {\rm d} s\\
  & \leq \frac{2}{\sqrt{3}} R_2^{1/2} |w(R_2,0)| + \left(\frac{2(R_2 - R)}{R}\right)^{1/2} \|G\|_{L^2(\{s: R<|s|<R_2\})}\\
  & \leq C \|\vec{w}(0)\|_{\mathcal{H}(R_2)} + \|G\|_{L^2(\{s: R<|s|<R_2\})}. 
 \end{align*}
 Here $C$ is an absolute constant. Combining this with Lemma \ref{lemma inequality}, we obtain 
  \begin{align*}
  \|\chi_{R} w\|_{Y(\Rm)} \leq C_1 \left( 2\|G\|_{L^2(\{s: R<|s|<R_2\})} + \|\chi_{R,R_2} w\|_{Y(\Rm)}^5+ \eta^4 \|\chi_{R,R_2} w\|_{Y(\Rm)} +C_2 \varepsilon \right).
 \end{align*}
Here $C_2> 1$ is an absolute constant. By choosing the same constants $\eta$, $\beta$ as in Lemma \ref{lemma connection} and applying a continuity argument in $R$, we obtain that if $\|G\|_{L^2(\{s: R<|s|<R_2\})} \leq \beta$ and $\varepsilon \leq \varepsilon_2 \doteq C_2^{-1} \beta$, then 
 \[
   \|\chi_{R,R_2} w\|_{Y(\Rm)} \leq \frac{8}{3}C_1 \left(2\|G\|_{L^2(\{s: R<|s|<R_2\})} +C_2\varepsilon \right) \leq 8C_1 \beta. 
 \]
 As in the proof of Lemma \ref{lemma connection}, an application of the nonlinear radiation profile shows
 \begin{align*}
 2\sqrt{2\pi} \|G\|_{L^2(\{s:|s|>R\})} & \leq  C_1\left(\|\chi_{R,R_2} w\|_{Y(\Rm)}^5 + \|\chi_{R,R_2} S\|_{Y(\Rm)}^4 \|\chi_{R,R_2} w\|_{Y(\Rm)}\right) + \varepsilon \\
 & \leq C_1 \left((8C_1 \beta)^4 + \eta^4\right) \|\chi_{R,R_2} w\|_{Y(\Rm)} + \varepsilon \nonumber \\
 & \leq \frac{3}{8C_1}\cdot \frac{8}{3}C_1 \left(2\|G\|_{L^2(\{s: R<|s|<R_2\})} + C_2 \varepsilon \right) + \varepsilon\\
 & \leq 2\|G\|_{L^2(\{s: R<|s|<R_2\})} +2 C_2 \varepsilon. 
\end{align*}
It immediately follows that 
\[
 \|G\|_{L^2(\{s: R<|s|<R_2\})} \leq \frac{2}{3} C_2 \varepsilon \leq \frac{2}{3}\beta. 
\]
A continuity argument then yields that $\|G\|_{L^2(\{s: R_1<|s|<R_2\})} \leq 2\beta/3$. As a result, the estimates given above hold for $R=R_1$. A combination of these estimates with Remark \ref{double tail G} finishes the proof. 
\end{proof}

\begin{corollary}\label{nc lemma}
 Given a positive integer $n$, assume that a radial solution $u$ satisfies case (b) of Proposition \ref{main tool}. Let $S = S_n$, $e = e_n(x,t)$ and $w = u - S$ introduced at the beginning of Section \ref{sec: pre}. For any $c\leq c_2$, there exists a constant $\delta_0=\delta_0(n,c)>0$, such that if $\delta < \delta_0$, then $u$ must be defined at least in the region $\Omega_{c\alpha_n^2}$ with 
 \[
  \|\vec{w}(\cdot,0)\|_{\mathcal{H}(c\alpha_n^2)} + \|\chi_{c\alpha_n^2} w\|_{Y(\Rm)} \lesssim_{n,c} \delta. 
 \]
\end{corollary}
\begin{proof}
The case $c=c_2$ has been verified. We start by choosing an absolute constant $\gamma\in (3/4,1)$ such that 
\[
 \sup_{R>0} \|\chi_{\gamma R,R}W\|_{Y(\Rm)} < \frac{\eta}{2n}. 
\]
It suffices to verify that if Corollary \ref{nc lemma} holds for a constant $c \leq c_2$, then it also holds for $\gamma c$. Our induction hypothesis implies that if $\delta < \delta_0(n,c)$ is sufficiently small, then a solution $u$ in case (b) and the associated solutions(functions) $S$, $w$, $e$  satisfy ($R_2 = c \alpha_n^2$)
\begin{align*}
 \|\vec{w}(0)\|_{\mathcal{H}(R_2)}  + \|\chi_{R_2}\left(F(u) - F(S)- e(x,t)\right)\|_{L^1 L^2(\Rm \times \Rm^3)} & \lesssim_{n,c} \delta; \\
  \|\chi_{0} e_n (x,t)\|_{L^1 L^2(\Rm \times \Rm^3)} & \lesssim_n \delta; \\
  \|\chi_{\gamma R_2, R_2} S\|_{Y(\Rm)} & \leq \frac{\eta}{2} + \delta. 
\end{align*}
Therefore if $\delta < \delta_0(n,\gamma c)$ is sufficiently small, we may apply Lemma \ref{lemma connection 2} to conclude that the following inequality holds for any radius $R_1$ with $R_1 >R'$ and $R_1 \geq \gamma R_2$. 
\begin{align*}
  \|\chi_{R_1} w\|_{Y(\Rm)} + \left\|\vec{w}(0)\right\|_{\mathcal{H}(R_1)} \lesssim_{n,\gamma c} \delta. 
 \end{align*}
 Here $R'$ is determined by the maximal domain $\Omega_{R'}$ of the asymptotically equivalent solution $u$. Please note that the right hand right does not depends on $R_1$. This implies that $\gamma R_2 > R'$, otherwise the uniform boundedness of $\|\chi_{R_1} w\|_{Y(\Rm)}$ (thus $\|\chi_{R_1} u\|_{Y(\Rm)}$) as $R_1 \rightarrow R'$ would give a contradiction. As a result, we may insert $R_1 = \gamma R_2 = \gamma c \alpha_n^2$ in the inequality above and finish the proof. 
\end{proof} 

In order to verify the inequality $|\alpha_{j+1}/\alpha_j| \lesssim_j \delta$, we need to refine the upper bound given above and prove the following lemma with explicit estimates. 

\begin{lemma} \label{lemma ratio}
 Let $u$, $S_n$, $w =w_n$ and $e=e_n$ be as in Corollary \ref{nc lemma}. If $\delta < \delta_0(n)$ is sufficiently small and the identity $r^{1/2} |w(r,0)| = \beta_1$ holds for some $r_\ast>R'$, then we must have $r_\ast \lesssim_n \delta^2 \alpha_n^2$. Here $\beta_1$ is the constant we chose in Step One. 
\end{lemma} 
\begin{proof}
 First of all, we let $c<c_2$ be a small constant, which would be determined later, apply Corollary \ref{nc lemma} and deduce that if $\delta < \delta_0(n,c)$, then $u$ is at least defined in $\Omega_{c\alpha_n^2}$ with
 \[
  \|\chi_{c\alpha_n^2} w\|_{Y(\Rm)} + \|\vec{w}(\cdot,0)\|_{\mathcal{H}(c\alpha_n^2)} + \|\chi_{c\alpha_n^2} (F(u)-F(w)-e_n)\|_{L^1 L^2} \lesssim_{c,n} \delta. 
 \]
 In addition, we have 
 \[
  \|\chi_0 e_n\|_{L^1 L^2} \lesssim_n \delta. 
 \]
 Without loss of generality, we may assume that $c$ is a sufficiently small number such that 
 \[
   \|\chi_{0,c} W\|_{Y(\Rm)} \leq \frac{\eta}{3n}.
 \]
 Now we assume that $C > C_1(n,c)$ is a large constant to be determined later, where the lower bound $C_1(n,c)$ guarantees that
 \[
  \|\chi_{c\alpha_n^2} w\|_{Y(\Rm)} + \|\vec{w}(\cdot,0)\|_{\mathcal{H}(c\alpha_n^2)} + \|\chi_{c\alpha_n^2} (F(u)-F(w)-e_n)\|_{L^1 L^2} \leq C \delta, 
 \]
 which implies 
 \[
  \sup_{r\geq c \alpha_n^2} r^{1/2} |w(r,0)| < C \delta;
 \]
 and define a sequence $c\alpha_n^2 > r_0 > r_1 > r_2 > \cdots$ by 
 \[
  r_k = \max\{r>0: r^{1/2} |w(r,0)| = 2^k C \delta\}, \qquad k=0,1,2, \cdots, \left\lfloor\log_2 \frac{\beta_1}{C \delta}\right\rfloor.
 \]
 When $\delta$ is sufficiently small, we may apply Lemma \ref{lemma connection} between radii $r_k$ and $c\alpha_n^2$ to deduce
 \[
  \|\vec{w}(0)\|_{\mathcal{H}(r_k)} + \|\chi_{r_k} w\|\lesssim_{n} 2^k C\delta, \qquad k\geq 0. 
 \]
 We recall Remark \ref{remark channel} and deduce 
 \begin{align*}
  \|G_\pm\|_{L^2(r_{k},r_{k-1})} & \lesssim_1 \left\|\chi_{r_{k}, r_{k-1}} \left(F(u) - F(S_n) - e_n\right)\right\|_{L^1 L^2} \\
  & \lesssim_n \delta + \left\|\chi_{r_{k}, r_{k-1}} \left(F(w+S_n) - F(S_n)\right)\right\|_{L^1 L^2}\\
  & \lesssim_n \delta + \left(\|\chi_{r_{k}, r_{k-1}} w\|_{Y(\Rm)}^4 + \|\chi_{r_{k}, r_{k-1}} S_n\|_{Y(\Rm)}^4\right) \|\chi_{r_{k}, r_{k-1}} w\|_{Y(\Rm)}\\
  & \lesssim_n \delta+ \left((2^k C \delta)^4 + \left(\frac{r_{k-1}}{\alpha_n^2}\right)^{2/5}\right) 2^k C \delta \\
  & \lesssim_n \left(C^{-1} 2^{-k} + (2^k C \delta)^4 + \left(\frac{r_{k-1}}{\alpha_n^2}\right)^{2/5}\right) 2^k C \delta. 
 \end{align*}
 Here $G_\pm$ is the corresponding radiation profile of $\vec{w}(0)$. By the definition of $r_k$ we have 
 \begin{align*}
  \left|r_k^{1/2} 2^k C \delta - r_{k-1}^{1/2} 2^{k-1} C\delta\right| & = \left|r_k |w(r_k,0)| - r_{k-1} |w(r_{k-1},0)|\right|\\
  & \leq  \left|r_k w(r_k,0) - r_{k-1} w(r_{k-1},0)\right|\\
  & \leq \int_{r_k}^{r_{k-1}} \left(|G_+(s)|+|G_-(s)|\right) {\rm d} s\\
  & \lesssim_n \left(C^{-1} 2^{-k} + (2^k C \delta)^4 + \left(\frac{r_{k-1}}{\alpha_n^2}\right)^{2/5}\right) 2^k C \delta (r_{k-1}-r_k)^{1/2}. 
 \end{align*}
 Thus 
 \begin{align*}
  \left|\frac{2r_{k}^{1/2}}{r_{k-1}^{1/2}} - 1\right|  \lesssim_n \left(C^{-1} 2^{-k} + (2^k C \delta)^4 + \left(\frac{r_{k-1}}{\alpha_n^2}\right)^{2/5}\right) \left(1-\frac{r_k^{1/2}}{r_{k-1}^{1/2}}\right).
 \end{align*}
 Therefore there exists a constant $c_n$ (depending on $n$ only) such that 
 \begin{align} \label{ratio estimate rk}
   \frac{2r_{k}^{1/2}}{r_{k-1}^{1/2}} \leq 1 + c_n\left( C^{-1} 2^{-k} + (2^k C \delta)^4 + \left(\frac{r_{k-1}}{\alpha_n^2}\right)^{2/5}\right).
 \end{align}
 Now we choose $c=c(n)$ and $C=C(n,c)=C(n)$ be sufficiently small/large constants, as well as an additional large constant $N=N(n) \in \mathbb{Z}^+$ such that 
 \begin{align*}
  &c_n c^{2/5} < \frac{1}{6};& & c_n/C < \frac{1}{6};& &c_n \beta_1^4 2^{-4N} < \frac{1}{6}. 
 \end{align*}
 As a result, if $\delta < \delta_0(n)$ is sufficiently small, then 
 \[
  \frac{2r_{k}^{1/2}}{r_{k-1}^{1/2}} < \frac{3}{2} \; \Longrightarrow \; \frac{r_{k}}{r_{k-1}} < \frac{9}{16}, \qquad \forall k = 1,2,\cdots,  K\doteq \left\lfloor\log_2 \frac{\beta_1}{C \delta}\right\rfloor - N. 
 \]
 Inserting these in \eqref{ratio estimate rk} yields ($k=1,2,\cdots,K$)
\begin{align*}
  \frac{2r_{k}^{1/2}}{r_{k-1}^{1/2}} & \leq 1 + \left( \frac{1}{6} \cdot 2^{-k} + c_n (2^k C \delta)^4 + \frac{1}{6}\left(\frac{9^{k-1}}{16^{k-1}}\right)^{2/5}\right)\\
 & \leq \exp  \left( \frac{1}{6} \cdot 2^{-k} + c_n (2^k C \delta)^4 + \frac{1}{6}\left(\frac{9^{k-1}}{16^{k-1}}\right)^{2/5}\right).
\end{align*}
It follows that 
\begin{align*}
 \frac{4^K r_K}{r_0} \leq \exp \sum_{k=1}^K  \left( \frac{1}{3} \cdot 2^{-k} + 2 c_n (2^k C \delta)^4 + \frac{1}{3}\left(\frac{9^{k-1}}{16^{k-1}}\right)^{2/5}\right) \lesssim_1 1. 
\end{align*}
 This finally gives
 \[
   r_\ast < r_K \leq 4^{-K} r_0 < 2^{2N+2} C^2 \beta_1^{-2} \delta^2 \cdot c \alpha_n^2 \lesssim_n \delta^2 \alpha_n^2
 \]
 and finish the proof. 
\end{proof}

\subsection{Step three}

In the last subsection we prove that if Proposition \ref{main tool} holds for a positive integer $n$, then it holds for $n+1$ as well. We start by choosing a small positive constant $c_1 = c_1(n)$ satisfying 
\[
 \|\chi_{0,c_1} W\|_{Y(\Rm)} \leq \frac{\eta}{3n}. 
\]
It suffices to consider the pairs $(v_L,u)$ with $\delta  = \|\chi_0 v_L\|_{Y(\Rm)}< \delta_0(n,c_1)$, where the upper bound is given by Corollary \ref{nc lemma} and solely depends on $n$. It is not difficult to see that we only need to consider solutions $u$ satisfying case (b) of the proposition for the positive integer $n$. In fact, if $u$ satisfies case (a) for the positive integer $n$, then it also satisfies case (a) for the positive integer $n+1$, with the same choice of $\alpha_j$'s. According to Corollary \ref{nc lemma} and the nonlinear estimates given in Remark \ref{nonlinear estimate}, the solution $u$ and associated solutions/functions $w_n$, $S_n$ and $e_n(x,t)$ defined at the beginning of Subsection \ref{sec: pre} satisfy 
\begin{align*}
 \left\|\vec{w}_n(0)\right\|_{\mathcal{H}(R_2)} + \|\chi_{R_2}\left(F(u) - F(S_n)- e_n(x,t)\right)\|_{L^1 L^2(\Rm \times \Rm^3)} & \lesssim_n \delta; \\
 \|\chi_0 e_n(x,t)\|_{L^1 L^2(\Rm \times \Rm^3)} & \lesssim_n \delta; \\
 \left\|\chi_{R_2} w_n\right\|_{Y(\Rm)} & \lesssim_n \delta.
\end{align*}
Here $R_2 = c_1 \alpha_n^2$ and $u$ must be defined in a maximal exterior region $\Omega_{R'}$ with $R' < R_2$. By further reducing the upper bound of $\delta$ if necessary, we see that the first inequality above also implies 
\begin{equation} \label{endpoint case R2 1}
 R_2^{1/2} |w_n(R_2,0)| \leq \|\vec{w}_n(0)\|_{\mathcal{H}(R_2)} < \frac{\beta_1}{10}.  
\end{equation}
Here $\beta_1$ is the absolute constant defined at the beginning of Step one. There are two cases:

\paragraph{Case one} In this case we assume 
\[
 \sup_{R'<r\leq R_2} r^{1/2} |w_n(r,0)| \leq \beta. 
\]
Our choice of $c_1$ implies that 
\[
 \|\chi_{0,R_2} S_n\|_{Y(\Rm)} \leq \frac{\eta}{3} + \delta. 
\]
A combination of the estimates given above implies that if $\delta < \delta_0(n+1)$ is sufficiently small, then we may apply Lemma \ref{lemma connection} for any interval $[R_1,R_2]$, as long as $R_1 > R'$, and obtain 
\[
  \|\chi_{R_1} w_n\|_{Y(\Rm)} + \left\|\vec{w}_n(0)\right\|_{\mathcal{H}(R_1)}  \lesssim_n R_1^{1/2} |w_n(R_1,0)| + \delta. 
\]
Again the uniform upper bound of the $Y$ norm for all $R_1>R'$ implies that $R'=0^-$. Therefore the estimate above also holds for $R_1 = 0$, which becomes 
\[
  \|\chi_{0} w_n\|_{Y(\Rm)} + \left\|\vec{w}_n(0)\right\|_{\dot{H}^1 \times L^2}  \lesssim_{n+1} \delta. 
\]
This implies that $u$ satisfies (a) for the positive integer $n+1$. 

\paragraph{Case two} In this case we have 
\[
 \sup_{R'<r\leq R_2} r^{1/2} |w_n(r,0)| > \beta. 
\]
Combining this with the continuity, the fact $\beta_1<\beta/2$ and \eqref{endpoint case R2 1}, we may always find a radius $R_1 \in (R',R_2)$ such that 
 \begin{align*}
  \sup_{R_1\leq r\leq R_2} r^{1/2} |w_n (r,0)| = R_1^{1/2} |w_n (R_1,0)| = \beta_1. 
 \end{align*}
 Now we choose $\alpha_{n+1} = \pm (c_2^{-1} R_1)^{1/2}$, where the sign is equal to that of $w_n (R_1,0)$. Clearly 
 \[
  |\alpha_{n+1}| = (c_2^{-1} R_1)^{1/2} < (c_2^{-1} R_2)^{1/2} = (c_1/c_2)^{1/2} |\alpha_n| < |\alpha_n|. 
 \]
 By Lemma \ref{lemma ratio}, when $\delta < \delta_0(n)$ is sufficiently small, we must have 
 \begin{align} \label{estimate of kappa nnplus}
  &R_1 \lesssim_n \delta^2 \alpha_n^2; & &\Longrightarrow& & \frac{|\alpha_{n+1}|}{|\alpha_n|} \lesssim_n \delta. 
 \end{align}
 Next we let 
 \[
  S_{n+1} =  v_L+ \sum_{j=1}^{n+1} W^{\alpha_j} = S_n + W^{\alpha_{n+1}}
 \]
 and define $e_{n+1}$ and $w_{n+1} = u - S_{n+1} = w_n - W^{\alpha_{n+1}}$ accordingly. Combining the estimates for $w_n$ and \eqref{estimate of kappa nnplus}, we observe that 
  \begin{align*}
  \|\vec{w}_{n+1}(0)\|_{\mathcal{H}(R_2)}  + \|\chi_{0} e_{n+1}(x,t)\|_{L^1 L^2(\Rm \times \Rm^3)} \qquad \qquad  & \\
 + \|\chi_{R_2}\left(F(u) - F(S_{n+1})- e_{n+1}(x,t)\right)\|_{L^1 L^2(\Rm \times \Rm^3)} & \lesssim_n \delta; \\
   \|\chi_{R_1,R_2} S_{n+1}\|_{Y(\Rm)} & \leq \frac{2\eta}{3} + \delta; \\
  \sup_{R_1\leq r\leq R_2} \left(r^{1/2} |w_{n+1}(r,0)|\right) & \leq 2\beta_1 < \beta;\\
  |w_{n+1}(R_1,0)| & = 0.
 \end{align*}
 As a result, if $\delta < \delta_0(n+1)$ is sufficiently small, then we may apply Lemma \ref{lemma connection} and obtain 
 \begin{align*}
   \|\chi_{R_1} w_{n+1}\|_{Y(\Rm)} + \left\|\vec{w}_{n+1}(\cdot,0)\right\|_{\mathcal{H}(R_1)} \lesssim_{n+1} \delta. 
 \end{align*}
This is the case (b) for positive integer $n+1$. 

\begin{remark}
 Given a positive integer $n$, we may determine the exact values of $J$ and parameters $\alpha_1, \alpha_2, \cdots, \alpha_J$ (in case (b) $J=n$) for any pair $(v_L, u)$ with a small norm $\|\chi_0 v_L\|_{Y(\Rm)} < \delta_0(n)$, by following the procedure given above. Please note that a small perturbation of $\alpha_j$'s may still satisfy the conditions given in Proposition \ref{main tool}. 
\end{remark}

\section{Proof of main theorem}

In this section we prove Theorem \ref{main thm}. We start by giving a lemma concerning free waves with highly concentrated radiation profiles. 

\begin{lemma} \label{calculation 1}
 Let $v_L$ be a radial free wave and $I = [a,b] \subset \Rm^+$ be an interval. Then 
 \begin{align*}
  \|\chi_0 v_L\|_{Y(\Rm)} + \|v_L\|_{Y(\Rm^+)} + \|v_L(\cdot,0)\|_{L^6(\Rm^3)} & \lesssim_1 \|G_+\|_{L^2(\Rm\setminus I)} + \left(\frac{b-a}{a}\right)^{1/2} \|G_+\|_{L^2(\Rm)}; \\
   \sup_{t>-a} \left\|(\nabla v_L (\cdot,t), \partial_t v_L (\cdot,t))\right\|_{L^2(\{x: |x|<a+t\})} &\lesssim_1 \|G_+\|_{L^2(\Rm\setminus I)}; \\
   \left\|\vec{v}_L (\cdot,0)\right\|_{\mathcal{H}(b)} & \lesssim_1 \|G_+\|_{L^2(\Rm\setminus I)} + \left(\frac{b-a}{a}\right)^{1/2} \|G_+\|_{L^2(\Rm)}; 
 \end{align*}
 Here $G_+$ is the radiation profile of $v_L$ in the positive time direction. 
 \end{lemma}
\begin{proof}
 The proof can be given by a straight-forward calculation. For convenience we let 
 \begin{align*}
  &\gamma = \frac{b-a}{a};& &\delta = \|G_+\|_{L^2(\Rm\setminus I)};& & M=\|G_+\|_{L^2(\Rm)}.
 \end{align*}
 Let $G_-$ be the radiation profile of $v_L$ in the negative time direction. We split $v_L$ into two parts:
  \[
   v_L = v_{L}^1 + v_L^2,
  \]
  whose radiation profiles $G_-^1$, $G_-^2$ in the negative time direction are given by 
  \begin{align*}
   &G_-^1(s) = \left\{\begin{array}{ll} G_-(s), & s\notin [-b,-a]; \\ 0, & s\in [-b, -a]. \end{array}\right.& 
   &G_-^2(s) = \left\{\begin{array}{ll} 0 , & s\notin [-b,-a]; \\ G_-(s), & s\in [-b, -a]. \end{array}\right.
  \end{align*}
  The symmetry \eqref{symmetry of G pm} implies that $\|G_-^1\|_{L^2(\Rm)} \leq \delta$ and $\|G_-^2\|_{L^2(\Rm)} \leq M$. By the Strichartz estimates we have 
   \begin{equation} \label{v1L Y norm}
   \|v_L^1\|_{Y(\Rm)} \lesssim_1 \delta.
  \end{equation}
  In order to estimate the norm $\|\chi_0 v_L^2\|_{Y(\Rm)}$, we recall the formula 
  \[
   v_L^2 (x,t) = \frac{1}{|x|} \int_{t-|x|}^{t+|x|} G_-^2 (s) {\rm d} s.  
  \]
  It follows that if $|x|<a+t$, then $v_L^2(x,t) = 0$; and that 
  \begin{equation} \label{pointwise vL2}
   |v_L^2 (x,t)| \leq \frac{1}{|x|} \int_{-b}^{-a} |G_-^2 (s)| {\rm d} s \lesssim_1 \frac{\gamma^{1/2} a^{1/2} M}{|x|}. 
  \end{equation} 
  Thus 
  \begin{align*}
   \|\chi_0 v_L^2\|_{Y(\Rm)}^5 & \lesssim_1 \int_{-\infty}^{-a/2} \left(\int_{|x|>|t|} \left|\frac{\gamma^{1/2} a^{1/2} M}{|x|}\right|^{10} {\rm d} x\right)^{1/2} {\rm d} t\\
    & \qquad + \int_{-a/2}^\infty \left(\int_{|x|>a+t} \left|\frac{\gamma^{1/2} a^{1/2} M}{|x|}\right|^{10}{\rm d} x\right)^{1/2} {\rm d} t \\
    & \lesssim_1 \gamma^{5/2} a^{5/2} M^5 \left(\int_{-\infty}^{-a/2} |t|^{-7/2} {\rm d} t + \int_{-a/2}^\infty (a+t)^{-7/2} {\rm d} t\right)\\
    & \lesssim_1 \gamma^{5/2} M^5. 
  \end{align*}
  In summary, we have 
  \begin{equation*}
    \|\chi_0 v_L\|_{Y(\Rm)} \lesssim_1 \delta + \gamma^{1/2} M. 
  \end{equation*} 
  Since $v_L^2(x,t) = 0$ for all $|x|<t$, we also have 
  \[
   \|v_L^2\|_{Y(\Rm^+)} \leq \|\chi_0 v_L^2\|_{Y(\Rm^+)}\lesssim_1 \gamma^{1/2} M; \qquad \Rightarrow \qquad \|v_L\|_{Y(\Rm^+)} \lesssim_1 \delta + \gamma^{1/2} M. 
  \] 
  A direct calculation also shows that 
  \begin{align*}
   \|v_L^2(\cdot,0)\|_{L^6}^6 \lesssim_1 \int_{|x|>a} \left(\frac{\gamma^{1/2} a^{1/2} M}{|x|}\right)^6 {\rm d} x \lesssim_1 \gamma^3 M^6.
  \end{align*}
  Thus 
  \[
   \|v_L(\cdot,0)\|_{L^6(\Rm^3)} \leq \|v_L^1 (\cdot,0)\|_{L^6(\Rm^3)} + \|v_L^2(\cdot,0)\|_{L^6(\Rm^3)} \lesssim_1 \delta + \gamma^{1/2} M. 
  \]
  We still need to verify the last two inequalities. Clearly we have 
\[
 \|\vec{v}_L^1(\cdot,0)\|_{\dot{H}^1\times L^2} \lesssim_1 \delta. 
\]
Now we consider $v_L^2$. We recall that $\vec{v}_L^2(x,t) = 0$ for $|x|<a+t$. The explicit formula given above also implies 
\[
 \vec{v}_L^2(x,0) = \left(\frac{1}{|x|} \int_{-b}^{-a} G_-^2 (s) {\rm d} s,0\right), \qquad |x|>b. 
\]
A straightforward calculation shows that
\[
 \|\vec{v}_L^2(\cdot,0)\|_{\mathcal{H}(b)} \lesssim_1 \gamma^{1/2} M. 
\] 
A combination of the estimates above finishes the proof. 
\end{proof}

\subsection{Soliton resolution with weak local radiation}

Next we incorporate Lemma \ref{calculation 1} into Proposition \ref{main tool}, utilize the continuity of energy norm and Lemma \ref{local boundedness} to show that the soliton resolution holds as long as the local radiation is sufficiently weak. We first introduce the notation for the asymptotic equivalent linear free wave. Let $u$ be a radial solution to (CP1) defined for all $t\geq 0$. In addition, we assume that $\|\vec{u}(0)\|_{\mathcal{H}(\bar{t})} \ll 1$ is sufficiently small for some $\bar{t}>0$. It follows that the radiation profiles $G_+ \in L^2(\Rm)$ and $G_-\in L^2([\bar{t},+\infty))$ are both well-defined for the non-linear solution $u$. We define $v_{\bar{t},L}$ to be the radial linear free wave with radiation profiles 
\begin{align*}
 & G_{\bar{t},+} (s) = G_+ (s-\bar{t}), \quad s>0; & & G_{\bar{t},-}(s) = G_-(s+\bar{t}), \quad s>0. 
\end{align*}
It is not difficult to see $v_{\bar{t},L}$ is asymptotically equivalent to the time-translated solution $u(x,t+\bar{t})$ in $\Omega_0$. We call $v_{\bar{t},L}$ the asymptotically equivalent linear free wave associated to $u$ at time $\bar{t}$. Next we prove

\begin{lemma} \label{pre soliton resolution}
 Let $E_0>E(W,0)$ be a constant. Then the following holds for sufficiently large constant $\ell > \ell(E_0)$ and sufficiently small constant $\delta_1 < \delta_1(E_0)$: If $u$ is a radial solution to (CP1) satisfying
 \begin{itemize}
  \item $u$ is defined for all $t\geq 0$ with a energy $E \in [E(W,0),E_0]$;
  \item The initial data $(u_0,u_1)$ and the radiation profile $G_+$ of $u$ satisfy 
  \begin{align*}
   &\|(u_0,u_1)\|_{\mathcal{H}(\ell^{-1}t_0)} \leq \delta_1;& &\|G_+\|_{L^2([-t_0,-\ell^{-1}t_0])} \leq \delta_1;
  \end{align*}
  here $t_0>0$ is a time;
 \end{itemize}
 then there exists a nonnegative integer $J$ and a sequence $(\zeta_j,\lambda_j) \in \{+1,-1\}\times \Rm^+$ for $j=1,2,\cdots,J$ satisfying
 \begin{align*}
  \frac{\lambda_1}{t_0}, \frac{\lambda_2}{\lambda_1}, \cdots, \frac{\lambda_J}{\lambda_{J-1}} \lesssim_{E_0} \left(\delta_1 + \ell^{-1/2} E_0^{1/2}\right)^2;   
 \end{align*}
    such that 
   \[
      \left\|\vec{u}(\cdot, t_0)-\sum_{j=1}^J \zeta_j (W_{\lambda_j},0) - \vec{v}_{t_0,L} (\cdot,0)\right\|_{\dot{H}^1\times L^2} \lesssim_{E_0} \delta_1 + \ell^{-1/2} E_0^{1/2}.
   \]
   Here $v_{t_0,L}$ is the asymptotically equivalent free wave of $u$ defined above. In addition, the energy $E$ and the norm $\|\vec{u}(t_0)\|_{\dot{H}^1\times L^2}$ satisfies 
   \begin{align*}
   & \left|\|\vec{u}(t_0)\|_{\dot{H}^1\times L^2}^2 - J \|W\|_{\dot{H}^1}^2 - 8\pi \|G_+\|_{L^2([-t_0,+\infty))}^2 \right| \\
    &\qquad + \left|E - J E(W,0) - 4\pi \|G_+\|_{L^2([-t_0,+\infty))}^2\right| \lesssim_{E_0}  \delta_1 + \ell^{-1/2} E_0^{1/2}.
   \end{align*}
\end{lemma}

\begin{proof}
 Assume that $(u,t_0)$ satisfies the assumptions in Lemma \ref{pre soliton resolution} for some large constant $\ell \gg 1$ and a small constant $\delta_1$. Without loss of generality we assume that $\delta_1 < \delta_0$, where $\delta_0$ is the small constant in Lemma \ref{local boundedness}. By Lemma \ref{local boundedness} and the exterior scattering, the following inequality holds for any $s_0\in \Rm^+$:
\[
 8\pi \|G_+\|_{L^2([-s_0,+\infty))}^2 = \lim_{t\rightarrow +\infty} \int_{|x|>t-s_0} |\nabla_{t,x} u(x,t)|^2 {\rm d} x \leq 6E + K_0.
\]
Thus 
\begin{equation} \label{upper bound Gplus}
 8\pi \|G_+\|_{L^2(\Rm)}^2 \leq 6 E + K_0\lesssim_1 E_0.
\end{equation}
In addition, the small data theory and finite speed of wave propagation implies that 
\begin{equation} \label{outer estimate ell} 
 \|\vec{u}(t)\|_{\mathcal{H}(|t|+\ell^{-1}t_0)} < 2\delta_1; \qquad \Longrightarrow \qquad \|G_\pm\|_{L^2([\ell^{-1}t_0,+\infty))} \lesssim_1 \delta_1. 
\end{equation}
For each $\bar{t} \in [t_0/3,t_0]$ we let $v_{\bar{t},L}$ be the corresponding asymptotically equivalent free wave. We apply Lemma \ref{calculation 1} on $v_{\bar{t},L}$ and $I = [\bar{t}-\ell^{-1} t_0, \bar{t}+\ell^{-1}t_0]$ to deduce
\begin{align*}
 \|\chi_0 v_{\bar{t},L}\|_{Y(\Rm)} &\lesssim_1 \|G_{\bar{t},+}\|_{L^2(\Rm\setminus I)} + \left(\frac{2\ell^{-1}t_0}{\bar{t} -\ell^{-1}t_0}\right)^{1/2} \|G_{\bar{t},+}\|_{L^2(\Rm)}\\
 & \lesssim_1 \|G_+\|_{L^2((-\bar{t},-\ell^{-1}t_0))} + \|G_+\|_{L^2([\ell^{-1}t_0,+\infty))} + \|G_-\|_{L^2([\bar{t},\infty))}+ \ell^{-1/2} E_0^{1/2}\\
 & \lesssim_1 \delta_1 + \ell^{-1/2} E_0^{1/2}. 
\end{align*}
For convenience we let $\delta = \delta_1 + \ell^{-1/2} E_0^{1/2}$.  Similarly we have 
\begin{align} 
 \|v_{\bar{t},L}(\cdot,0)\|_{L^6(\Rm^3)} & \lesssim_1 \delta; \label{L6 estimate ell}\\
  \left\|(\nabla v_{\bar{t},L} (\cdot,0), \partial_t v_{\bar{t},L} (\cdot,0))\right\|_{L^2(\{x: |x|<\bar{t}-\ell^{-1}t_0\;\hbox{or}\; |x|>\bar{t}+\ell^{-1}t_0\})} & \lesssim_1 \delta. \label{energy estimate ell}
\end{align}
Now we choose a positive integer $n = n(E_0)$ satisfying 
\[
 (n-1)\|W\|_{\dot{H}^1}^2 > 6E_0 + K_0 + 1,
\]
then apply Lemma \ref{main tool} on the time-translated solution $u(x,\bar{t}+t)$ and asymptotically equivalent free wave $v_{\bar{t},L}$. It follows from the upper bound of $\|\chi_0 v_{\bar{t},L}\|_{Y(\Rm)}$ given above that if $\ell > \ell(E_0)$ is sufficiently large and $\delta_1<\delta_1(E_0)$ is sufficiently small, then for each $\bar{t}$ either of the following holds:
\begin{itemize}
 \item[(a)] There exists a sequence $\alpha_j(\bar{t})$ for $j=1,2,\cdots,J(\bar{t})$ with $0\leq J(\bar{t})<n$ satisfying
 \begin{align} \label{ratio estimate a}
 |\alpha_{j+1}(\bar{t})|/|\alpha_j(\bar{t})| &\lesssim_{E_0} \delta, \quad j =1,2,\cdots,J(\bar{t})-1; 
 \end{align}
 such that 
   \begin{equation} \label{norm estimate case an}
     \left\|\vec{u}(\cdot, \bar{t})-\sum_{j=1}^{J(\bar{t})} (W^{\alpha_j(\bar{t})},0) - \vec{v}_{\bar{t},L} (\cdot,0)\right\|_{\dot{H}^1\times L^2} \lesssim_{E_0} \delta.
   \end{equation}
 \item[(b)] There exists a sequence $\alpha_j(\bar{t})$ for $j=1,2,\cdots,n$ satisfying
 \begin{align} \label{ratio estimate b}
 |\alpha_{j+1}(\bar{t})|/|\alpha_j(\bar{t})| &\lesssim_{E_0} \delta, \quad j =1,2,\cdots,n-1;
 \end{align}
 such that 
   \[
     \left\|\vec{u}(\cdot, \bar{t})-\sum_{j=1}^{n} (W^{\alpha_j(\bar{t})},0) - \vec{v}_{\bar{t},L} (\cdot,0)\right\|_{\mathcal{H}(c_2 \alpha_n(\bar{t})^2)} \lesssim_{E_0} \delta.
   \]
\end{itemize}
We claim that unless $(u,\bar{t})$ satisfies case (a) with $J=0$, the following inequality also holds: 
\[
 \frac{\alpha_1 (\bar{t})^2}{\bar{t}} \lesssim_{E_0} \delta^2.
\]
In fact, a combination of \eqref{outer estimate ell}, \eqref{energy estimate ell} and the soliton resolution implies that 
\begin{align} 
  \left\|\sum_{j=1}^{J(\bar{t})} (W^{\alpha_j(\bar{t})},0) \right\|_{\mathcal{H}\left(\bar{t}+\ell^{-1} t_0\right)} & \lesssim_{E_0} \delta, & & \hbox{Case (a);} \label{upper bound sum W J1}\\
 \left\|\sum_{j=1}^{n} (W^{\alpha_j(\bar{t})},0) \right\|_{\mathcal{H}\left(\max\{\bar{t}+\ell^{-1} t_0, c_2 \alpha_n(\bar{t})^2\}\right)} & \lesssim_{E_0} \delta, & & \hbox{Case (b).} \label{upper bound sum W J2}
\end{align}
Combining the ratio estimate \eqref{ratio estimate a} or \eqref{ratio estimate b} with the fact 
\[
 \|(W^{\alpha},0)\|_{\mathcal{H}(c\alpha^2)} \simeq_1 c^{-1/2}, \qquad c\geq 1,
\]
we obtain the following estimate when $\delta_1 \leq \delta_1(E_0)$ is small and $\ell > \ell(E_0)$ is large: 
\begin{equation} \label{lower bound sum W J}
 \left\|\sum_{j=1}^{J(\bar{t})} (W^{\alpha_j(\bar{t})}, 0)\right\|_{\mathcal{H}(c \alpha_1(\bar{t})^2)}  \gtrsim_1 c^{-1/2}, \qquad c\geq 1. 
\end{equation} 
A comparison of \eqref{upper bound sum W J1} (or \eqref{upper bound sum W J2}) and \eqref{lower bound sum W J} immediately verifies our claim. Next we claim that the following energy estimates hold for each time $\bar{t}$ as well, as long as  $\ell_1 > \ell_1(E_0)$ is sufficiently large and $\delta_1<\delta_1(E_0)$ is sufficiently small.
 \begin{align*}
    \left|\|\vec{u}(\bar{t})\|_{\dot{H}^1\times L^2}^2 - J(\bar{t}) \|W\|_{\dot{H}^1}^2 - 8\pi \|G_+\|_{L^2([-\bar{t},+\infty))}^2 \right| & \lesssim_{E_0} \delta; & & \hbox{(Case a)}\\
    \left|E - J(\bar{t}) E(W,0) - 4\pi \|G_+\|_{L^2([-\bar{t},+\infty))}^2\right| & \lesssim_{E_0} \delta; & & \hbox{(Case a)} \\
    \|\vec{u}(\bar{t})\|_{\dot{H}^1\times L^2}^2 & \geq (n-1) \|W\|_{\dot{H}^1}^2. & & \hbox{(Case b)}
   \end{align*}
In fact, the scaling separation given by the argument above and the energy localization \eqref{energy estimate ell} of $v_{\bar{t},L}$ yields the followsing almost orthogonality estimates: 
\begin{align*}
 \int_{\Rm^3} \left|\nabla W^{\alpha_j(\bar{t})}(x) \cdot \nabla W^{\alpha_k(\bar{t})}(x)\right| {\rm d} x  & \lesssim_{E_0} \delta, \qquad j\neq k; \\
 \int_{\Rm^3} \left|\nabla W^{\alpha_j(\bar{t})}(x) \cdot \nabla v_{\bar{t},L} (x,0)\right| {\rm d} x & \lesssim_{E_0} \delta.
\end{align*}
In addition, we have 
\begin{align*}
 0\leq \|\vec{v}_{\bar{t},L} (\cdot,0)\|_{\dot{H}^1\times L^2}^2 - 8\pi \|G_+\|_{L^2([-\bar{t},+\infty))}^2  & = 8\pi \|G_-\|_{L^2([\bar{t},+\infty))}^2  \lesssim_1 \delta_1^2. 
\end{align*}
A combination of these estimates and \eqref{norm estimate case an} shows that in case (a) we have 
\begin{align}
& \left|\|\vec{u}(\bar{t})\|_{\dot{H}^1\times L^2}^2 - J(\bar{t}) \|W\|_{\dot{H}^1}^2 - 8\pi\|G_+\|_{L^2([-\bar{t},+\infty))}^2 \right| \nonumber\\
 & \leq  \left|\|\vec{u}(\bar{t})\|_{\dot{H}^1\times L^2}^2 \!-\! \sum_{j=1}^{J(\bar{t})} \|(W^{\alpha_j(\bar{t})},0)\|_{\dot{H}^1\times L^2}^2 - \|\vec{v}_{\bar{t},L} (\cdot,0)\|_{\dot{H}^1\times L^2}^2 \right| \nonumber\\
 &\qquad + \left|\|\vec{v}_{\bar{t},L} (\cdot,0)\|_{\dot{H}^1\times L^2}^2 - 8\pi \|G_+\|_{L^2([-\bar{t},+\infty))}^2 \right| \lesssim_{E_0} \delta. \label{norm estimate case aaa} 
\end{align}
Similarly in case (b) we have 
\begin{align*}
 & \left|\|\vec{u}(\bar{t})\|_{\mathcal{H}(c_2 \alpha_n(\bar{t})^2)}^2 - (n-1) \|W\|_{\dot{H}^1}^2 - \|\nabla W\|_{L^2(\{x: |x|>c_2\})}^2 - \|\vec{v}_{\bar{t},L} (0)\|_{\mathcal{H}(c_2 \alpha_n(\bar{t})^2)}^2 \right| \\
 & \leq  \left|\|\vec{u}(\bar{t})\|_{\mathcal{H}(c_2 \alpha_n(\bar{t})^2)}^2 - \sum_{j=1}^{n} \left\|\nabla W^{\alpha_j(\bar{t})}\right\|_{L^2(\{x: |x|>c_2 \alpha_n(\bar{t})^2\})}^2  - \|\vec{v}_{\bar{t},L} (0)\|_{\mathcal{H}(c_2 \alpha_n(\bar{t})^2)}^2 \right| \\
 & \qquad +  \sum_{j=1}^{n-1} \left\|\nabla W^{\alpha_j(\bar{t})}\right\|_{L^2(\{x: |x|<c_2 \alpha_n(\bar{t})^2\})}^2 \lesssim_{E_0} \delta. 
\end{align*}
This immediately gives $\|\vec{u}(\bar{t})\|_{\dot{H}^1\times L^2}^2 \geq (n-1) \|W\|_{\dot{H}^1}^2$ for sufficiently large $\ell > \ell(E_0)$ and sufficiently small $\delta_1 < \delta_1(E_0)$. Finally we consider the energy of the nonlinear wave equation (in case (a) only)
\[
 E(u,u_t) = \int_{\Rm^3} \left(\frac{1}{2} |\nabla u|^2 + \frac{1}{2} |u_t|^2 - \frac{1}{6} |u|^6 \right) {\rm d} x = \frac{1}{2} \|\vec{u}\|_{\dot{H}^1\times L^2}^2 - \frac{1}{6} \int_{\Rm^3} |u|^6 {\rm d} x 
\]
and prove that
\[
   \left|E - J(\bar{t}) E(W,0) - 4\pi \|G_+\|_{L^2([-\bar{t},+\infty))}^2\right| \lesssim_{E_0} \delta.
\]
In view of \eqref{norm estimate case aaa}, it suffices to show 
\begin{equation} \label{to prove L6}
 \left|\int_{\Rm^3} |u(x,\bar{t})|^6 {\rm d} x - J(\bar{t}) \int_{\Rm^3} |W(x)|^6 {\rm d} x\right| \lesssim_{E_0} \delta. 
\end{equation}
We write
\[
 u(x,\bar{t}) = \sum_{j=1}^{J(\bar{t})} W^{\alpha_j(\bar{t})} (x) + v_{\bar{t},L} (x,0) + w(x,\bar{t})
\]
and observe  
\begin{itemize}
 \item $\|w(x,\bar{t})\|_{L^6 (\Rm^3)} \lesssim_1 \|w(\cdot,\bar{t})\|_{\dot{H}^1(\Rm^3)} \lesssim_{E_0} \delta$ by \eqref{norm estimate case an}; 
 \item $\|v_{\bar{t},L} (\cdot,0)\|_{L^6} \lesssim_{1} \delta$ by \eqref{L6 estimate ell};  
 \item $\|W^{\alpha}\|_{L^6(\Rm^3)}$ is independent of $\alpha\neq 0$. In addition, scale separation implies that 
 \[
   \int_{\Rm^3} \left|W^{\alpha_j(\bar{t})}\right|^{m} \left|W^{\alpha_k(\bar{t})}\right|^{6-m} {\rm d} x \lesssim_{E_0} \delta, \qquad j\neq k, \; m = 1,2,\cdots,5. 
 \]
 The estimates for interaction terms with three or more scales immediately follow from an interpolation. 
\end{itemize}
A combination of these estimate then verifies our claim on energy estimates. In summary, if $\ell > \ell(E_0)$ is sufficiently large, $\delta_1 < \delta_1(E_0)$ is sufficiently small, $u$ and $t_0$ satisfy the conditions of Lemma \ref{pre soliton resolution}, then for each $\bar{t} \in [t_0/3,t_0]$ either of the following holds:
\begin{itemize}
 \item [(a1)] The conclusion of Lemma \ref{pre soliton resolution} holds at time $\bar{t}$ with $\zeta_j = {\rm sign} (\alpha_j(\bar{t}))$ and $\lambda_j(\bar{t}) = \alpha_j(\bar{t})^2$; 
 \item[(b1)] The energy norm satisfies $\|\vec{u}(\bar{t})\|_{\dot{H}^1\times L^2}^2 \geq (n-1) \|W\|_{\dot{H}^1}^2 > 6 E_0 + K_0 + 1$. 
\end{itemize}
By Lemma \ref{local boundedness}, the case (a1) holds for at least one time $\bar{t}\in [t_0/3,t_0]$. In addition, in view of the identity $\|W\|_{\dot{H}^1}^2 = 3 E(W,0)$, the energy estimate above also implies that the following inequality holds in case (a1):
\[
 \|\vec{u}(\bar{t})\|_{\dot{H}^1\times L^2}^2 \leq 3E + C(E_0) \delta. 
\]
A continuity argument then verifies that case (a1) holds for all time $\bar{t} \in [t_0/3,t_0]$, as long as $\delta_1<\delta_1(E_0)$ is sufficiently small and $\ell>\ell(E_0)$ is sufficiently large. Finally we choose $J=J(t_0)$, $\lambda_j = |\alpha_j(t_0)|^2$, $\zeta_j = \hbox{sign} (\alpha_j(t_0))$ and finish the proof. 
\end{proof}
 
\subsection{Proof of Theorem \ref{main thm}} 

Now we are ready to give the proof of the main theorem. Let $\delta$ be a small constant, which will be determined later, and $\ell = E_0 \delta^{-2} \ll 1$ be a large constant accordingly. Now let $u$ be a solution as in the main theorem and $R$ be a large radius such that 
\[
 \|\vec{u}(0)\|_{\mathcal{H}(R)} < \delta.
\]
We start by considering the local radiation strength function $\varphi_\ell (t) \doteq \|G_+\|_{L^2([-t, -\ell^{-1} t])}$ for $t \geq \ell R$. We define the set of time with weak local radiation strength:
\[
 Q = \{t>\ell R: \varphi_\ell (t) < \delta\}. 
\]
By the continuity of $\varphi_\ell$ and the fact $G_+\in L^2(\Rm)$, the set $Q$ is an open subset of $\Rm$ containing a neighbourhood of $+\infty$. We first verify that the soliton resolution holds for $t \in \bar{Q}$, the closure of $Q$. Assume that $t \in \bar{Q}$ (thus $\varphi_\ell(t) \leq \delta$). An application of Lemma \ref{pre soliton resolution} immediately gives a sequence $\{(\zeta_j(t),\lambda_j(t))\}_{j=1,2,\cdots,J(t)}$ with
   \begin{align*}
    & \max\left\{\frac{\lambda_1(t)}{t}, \frac{\lambda_2(t)}{\lambda_1(t)}, \cdots, \frac{\lambda_{J(t)}(t)}{\lambda_{J(t)-1}(t)}\right\} \lesssim_{E_0} \delta^2;&  & \zeta_j(t) \in \{+1,-1\};
   \end{align*}
   such that 
   \begin{equation} \label{soliton resolution 12}
     \left\|\vec{u}(\cdot, t)-\sum_{j=1}^{J(t)} \zeta_j(t) (W_{\lambda_j(t)},0) - \vec{v}_{t,L} (\cdot,0)\right\|_{\dot{H}^1\times L^2} \lesssim_{E_0} \delta;
   \end{equation}
  as long as $\delta < \delta(E_0)$ is sufficiently small. In addition, the energy $E$ of $u$ satisfies 
   \begin{align} \label{energy estimate 12}
    \left|E - J(t) E(W,0) - 4\pi \|G_{+}\|_{L^2([-t,+\infty))}^2\right| \lesssim_{E_0} \delta.
   \end{align}
This energy estimate implies that $J(t)$ is a non-increasing function of $t\in \bar{Q}$ and 
 \begin{equation} \label{upper bound of Jt}
  J(t) \leq \frac{E + C(E_0) \delta}{E(W,0)}. 
 \end{equation}
 \paragraph{Stable time periods} Let $J_1 > J_2 > \cdots > J_m$ be all possible values of $J(t)$ for $t\in Q$. We may split $Q$ into a few parts 
 \begin{align*}
  & Q = \bigcup_{k=1}^m Q_k; & & Q_k = \{t\in Q: J(t) = J_k\}. 
 \end{align*}
 By the non-increasing property of $J(t)$, the inequality $t_1<t_2$ holds if $t_1 \in Q_{k_1}$, $t_2 \in Q_{k_2}$ and $k_1<k_2$. It is not difficult to see that $Q_k$ are all nonempty open sets. We choose 
 \begin{align*}
  &a_k = \inf Q_k \geq \ell R;& &b_k = \sup Q_k > a_k.
 \end{align*}
 Clearly we have $\ell R \leq a_1 < b_1 < a_2 < b_2 < \cdots < a_m < b_m = +\infty$. Since $a_k \in \bar{Q}$, the approximated soliton resolution described above holds at time $a_k$. By continuity, the bubble number at time $a_k$ satisfies $J(a_k) = J_k$. The same happens at time $b_k$ for $k=1,2,\cdots, m-1$. According to \eqref{energy estimate 12}, the following $L^2$ estimate holds 
 \begin{equation} \label{conclusion 3 add}
  \|G_+\|_{L^2((-b_k,-a_k])} \lesssim_{E_0} \delta^{1/2}. 
 \end{equation}
 As a result, for any time $t\in [a_k,b_k]$ (please note that the right hand endpoint $+\infty$ is removed for $k=m$ here and in the argument below), we have 
 \[
  \|G_+\|_{L^2([-t,-\ell^{-1} t])} \leq \|G_+\|_{L^2((-b_k,-a_k])} + \|G_+\|_{L^2([-a_k,-\ell^{-1} a_k])} \lesssim_{E_0} \delta^{1/2}. 
 \]
 When $\delta<\delta(E_0)$ is sufficiently small, we may apply Lemma \ref{pre soliton resolution} again and conclude that there exists a sequence $\{(\zeta_j(t),\lambda_j(t))\}_{j=1,2,\cdots,J_k}$ for each time $t\in [a_k,b_k]$ with
   \begin{align} \label{conclusion 1} 
    & \max\left\{\frac{\lambda_1(t)}{t}, \frac{\lambda_2(t)}{\lambda_1(t)}, \cdots, \frac{\lambda_{J_k}(t)}{\lambda_{J_k-1}(t)}\right\} \lesssim_{E_0} \delta;&  & \zeta_j(t) \in \{+1,-1\};
   \end{align}
   such that 
   \begin{equation*} 
     \left\|\vec{u}(\cdot, t)-\sum_{j=1}^{J_k} \zeta_j(t) (W_{\lambda_j(t)},0) - \vec{v}_{t,L} (\cdot,0)\right\|_{\dot{H}^1\times L^2} \lesssim_{E_0} \delta^{1/2};
   \end{equation*}
Thus for each $k$ and $1\leq j \leq J_k$, we may choose 
 \begin{align}
  &\lambda_{k,j} (t) = \lambda_j(t), \quad t\in [a_k,b_k]; & & \zeta_{k,j} = \zeta_{j}(t).
 \end{align}
 Here the bubble number and signs $\zeta_{j}(t)\in \{\pm 1\}$ are independent of $t\in [a_k,b_k]$ by continuity. We still need to substitute $\vec{v}_{t,L}$ by a linear free wave independent of $t$ for each stable time period. We let $v_{k,L}$ be the linear free wave with the following radiation profile in the positive time direction:
 \[
  G_{k,+} (s) = \left\{\begin{array}{ll} G_+(s), & s> -b_k; \\ 0, & s < -b_k. \end{array}\right.
 \]
 Thus the time-translated version $v_{k,L} (x, \cdot + t)$ comes with a radiation profile 
 \[
  G_{k, t, +} (s) = \left\{\begin{array}{ll} G_+(s-t), & s> -b_k+t; \\ 0, & s < -b_k + t. \end{array}\right.
 \]
 Comparing the radiation profiles, we have 
 \begin{align*}
  \|\vec{v}_{t,L}(\cdot,0)-\vec{v}_{k,L}(\cdot,t)\|_{\dot{H}^1 \times L^2}^2 & = 8\pi \|G_{k,t,+} - G_{t,+}\|_{L^2(\Rm)}^2  = 8\pi \|G_{k,t,+} - G_{t,+}\|_{L^2((-\infty,0])}^2 \\
  &\lesssim_1 \|G_+\|_{L^2((-b_k,-t])}^2 + \|G_-\|_{L^2([t,+\infty))}^2 \lesssim_{E_0} \delta. 
 \end{align*}
As a result, we have 
 \begin{equation} \label{conclusion 2}
    \left\|\vec{u}(\cdot, t)-\sum_{j=1}^{J_k} \zeta_{k,j} (W_{\lambda_{k,j}(t)},0) - \vec{v}_{k,L} (\cdot,t)\right\|_{\dot{H}^1\times L^2} \lesssim_{E_0} \delta^{1/2}, \quad t\in [a_k,b_k]. 
 \end{equation}
In addition, we apply Lemma \ref{calculation 1} on $v_{k,L}(\cdot,\cdot+a_k)$ with $I = [a_k-\ell^{-1}a_k,a_k+\ell^{-1}a_k]$ to deduce 
  \begin{align}
  \|\chi_{|x|>|t-a_k|} v_{k,L}\|_{Y(\Rm)} & + \|v_{k,L}\|_{Y([a_k,+\infty))} + \sup_{t>\ell^{-1} a_k} \|(\nabla v_{k,L}, \partial_t v_{k,L})\|_{L^2(\{x: |x|<t-\ell^{-1} a_k\})}  \nonumber\\
  & \lesssim_1 \ell^{-1/2} E_0^{1/2} + \|G_+\|_{L^2((-b_k,-\ell^{-1}a_k]\cup[\ell^{-1}a_k,+\infty))} \nonumber \\
  & \lesssim_{E_0} \delta^{1/2}. \label{conclusion 3} 
 \end{align}
\paragraph{Collision periods} Now let us consider the collision time periods $[b_k, a_{k+1}]$. By the choice of $a_k$, $b_k$ and the continuity of $\varphi_\ell$, we must have 
 \[
  \varphi_\ell (b_k) = \varphi_\ell (a_{k+1}) = \delta, \qquad k =1, 2, \cdots, m-1. 
 \]
 The way we choose $[a_k,b_k]$ guarantees that $\varphi_\ell (t) \geq \delta$ for $t\in [b_k,a_{k+1}]$. Now let us give an upper bound of the ratio $a_{k+1}/b_k$. First of all, the fact $b_k, a_{k+1} \in \bar{Q}$ and the energy estimate \eqref{energy estimate 12} gives 
 \begin{align*}
   \left| 4\pi \int_{-a_{k+1}}^{-b_k} |G_+(s)|^2 {\rm d} s - p_k E(W,0)\right| \lesssim_{E_0}  \delta; & & p_k = J_k - J_{k+1} \in \mathbb{N}. 
 \end{align*}
 We may combine this upper bound with \eqref{upper bound of Jt} and the lower bound $\varphi_\ell(t) \geq \delta$ to deduce 
 \begin{equation} \label{conclusion 4}
  \delta^2 \left\lfloor \log_\ell \frac{a_{k+1}}{b_k} \right\rfloor \leq \int_{-a_{k+1}}^{-b_k} |G_+(s)|^2 {\rm d} s \leq \frac{E + C(E_0) \delta}{4\pi},
 \end{equation} 
as long as $\delta < \delta(E_0)$ is sufficiently small.  This helps us give an upper bound of the length for collision periods. Similarly we may give the upper bound of the ratio $a_{1}/R$. In fact, if $a_1 > \ell R$, then the way we choose $[a_k,b_k]$ guarantees that 
 \[
  \varphi_\ell (t) \geq \delta, \qquad t\in [\ell R, a_1]. 
 \]
 We may combine this with \eqref{energy estimate 12} to deduce 
 \begin{equation} \label{conclusion 5}
  \delta^2\left\lfloor \log_\ell \frac{a_1}{R} \right\rfloor \leq \int_{-a_1}^{-R} |G_+(s)|^2 {\rm d} s \leq \int_{-a_1}^{\infty} |G_+(s)|^2 {\rm d} s \leq \frac{E + C(E_0)\delta}{4\pi}.
 \end{equation} 
\paragraph{Completion of proof} Finally we may choose a sufficiently small parameter $\delta = \delta(E_0,\varepsilon)$ and complete the proof. The soliton resolution properties in stable time intervals follows from \eqref{conclusion 1}, \eqref{conclusion 2} and \eqref{conclusion 3}. The upper bounds for $a_{k+1}/b_k$ and $a_1/R$ follows from \eqref{conclusion 4}, \eqref{conclusion 5}. The radiation strength estimates for stable/collision time intervals follows from \eqref{energy estimate 12} and \eqref{conclusion 3 add}. In particular, if $\varepsilon < \varepsilon(E_0)$ is sufficiently small, we may pick 
\begin{align*}
 &\delta = c \varepsilon^2; & & \ell = E_0 c^{-2} \varepsilon^{-4}; & & L = \left(E_0 c^{-2} \varepsilon^{-4}\right)^\frac{E_0}{2\pi c^2 \varepsilon^4}.  
\end{align*}
Here $c= c(E_0)$ is a constant depending on $E_0$ only. 
 
 \section{One-pass theorem of multi-bubble solutions}
 
In this section we give a one-pass theorem of pure multi-bubble solutions. The first one-pass type theorem for the nonlinear wave equation was introduced by Grieger-Nakanishi-Schalg \cite{aboveground1}, as far as the author knows. Their theorem discussed the dynamics of solutions near the ground states, while the theory given in this section applies to solutions near pure multi-bubble solutions. We first introduce a few definitions.  Given a positive integer $n$, two small constants $\varepsilon, \kappa$, we define a ``neighbourhood of pure $k$-bubble'' $\mathcal{M}_n(\varepsilon, \kappa)$ to be the following subset of $\mathcal{H}$
 \[
  \mathcal{M}_n (\varepsilon,\kappa) = \left\{(u_0,u_1)\in \mathcal{H} : \begin{array}{c} \hbox{There exist}\; (\zeta_1, \lambda_1), \cdots, (\zeta_n,\lambda_n) \in \{\pm 1\} \times \Rm^+,  (w_0,w_1)\in \mathcal{H}, \\ 
  \hbox{with}\;  \lambda_{j+1}/\lambda_{j} < \kappa^2, j=1,\cdots, n-1; \; \|(w_0,w_1)\|_{\mathcal{H}} < \varepsilon;\\ 
  \hbox{such that}\; \displaystyle (u_0,u_1) = \sum_{j=1}^n \zeta_j (W_{\lambda_j},0) + (w_0,w_1).\end{array} \right\}.
 \]
 Clearly $\mathcal{M}_n(\varepsilon, \kappa)$ is an open subset of $\mathcal{H}$. In particular, we may fix $\{\zeta_j\}_{1\leq j\leq n} \in \{+1,-1\}^n$ and define 
 \[
  \mathcal{M}_n (\varepsilon,\kappa, \{\zeta_j\}_j) = \left\{(u_0,u_1)\in \mathcal{H} : \begin{array}{c} \hbox{There exist}\; \lambda_1, \cdots, \lambda_n \in \Rm^+,  (w_0,w_1)\in \mathcal{H}, \\ 
  \hbox{with}\;  \lambda_{j+1}/\lambda_{j} < \kappa^2, j=1,\cdots, n-1; \; \|(w_0,w_1)\|_{\mathcal{H}} < \varepsilon;\\ 
  \hbox{such that}\; \displaystyle (u_0,u_1) = \sum_{j=1}^n \zeta_j (W_{\lambda_j},0) + (w_0,w_1).\end{array} \right\}.
 \]
Clearly we have 
\[
 \mathcal{M}_n (\varepsilon,\kappa) = \bigcup_{\{\zeta_j\}_j \in \{+1,-1\}^n}  \mathcal{M}_n (\varepsilon,\kappa, \{\zeta_j\}_j). 
\]
In addition, if $\kappa < \kappa_0$ and $\varepsilon < \varepsilon_0$ are sufficiently small, then the sets $ \mathcal{M}_n (\varepsilon,\kappa, \{\zeta_j\}_j)$ and  $\mathcal{M}_n (\varepsilon,\kappa, \{\zeta'_j\}_j)$ are disjoint unless $\zeta_j = \zeta'_j$ for all $1\leq j \leq n$. Please note that the numbers $\kappa_0$ and $\varepsilon_0$ do not depend on $n$. 

 Now we give another way to define a roughly equivalent neighbourhood of pure $k$-bubble by considering the radiation. Given a sufficiently small constant $\delta > 0$, we let $\mathcal{R}(\delta)$ be the set of all radial initial data $(u_0,u_1)\in \dot{H}^1\times L^2(\Rm^3)$ such that the corresponding exterior solution $u$ defined in the exterior region $\Omega_0$ satisfies the following 
 \begin{itemize}
  \item The exterior solution $u$ is defined for all time $t\in \Rm$ such that $\|\chi_0 u\|_{Y(\Rm)} < +\infty$; 
  \item The nonlinear radiation profiles $G_\pm (s)$ satisfy $\|G_\pm\|_{L^2(\Rm^+)} < \delta$. 
 \end{itemize} 
 Lemma \ref{perturbation lemma} guarantees that $\mathcal{R}(\delta)$ is an open subset of $\mathcal{H}$. Next we define 
 \[
  \mathcal{R}_n(\delta) = \left\{(u_0,u_1)\in \mathcal{R}(\delta): (n-1/2) \|W\|_{\dot{H}^1}^2 < \|(u_0,u_1)\|_{\dot{H}^1\times L^2}^2<(n+1/2) \|W\|_{\dot{H}^1}\right\}
 \]
 for any positive integer $n$. It is clear that $\mathcal{R}_{n} (\delta)$ is also an open subset of $\mathcal{H}$. To see why these two kinds of neighbourhood are roughly equivalent, we need to apply the following lemma. 
 
 \begin{lemma} \label{two equivalent neig}
  Let $n$ be a fixed positive integer. Then given any $\delta > 0$, there exist $\kappa=\kappa(n,\delta)$ and $\varepsilon= \varepsilon(n, \delta)$ such that $\mathcal{M}_{n}(\kappa, \varepsilon) \subseteq \mathcal{R}_n (\delta)$. Conversely, given any $\kappa, \varepsilon> 0$, there exists $\delta = \delta(n,\kappa, \varepsilon)$, such that $\mathcal{R}_n (\delta) \subseteq \mathcal{M}_{n}(\kappa, \varepsilon)$. 
 \end{lemma}
 \begin{proof}
  Let us first assume that $(u_0,u_1) \in \mathcal{M}_n(\kappa, \varepsilon)$ can be given in the following form 
  \[
   (u_0,u_1) = \sum_{j=1}^n \zeta_j (W_{\lambda_j},0) + (w_0,w_1)
  \]
  with 
  \begin{align*}
   &\frac{\lambda_{j+1}}{\lambda_{j}} < \kappa^2, \; j=1,2,\cdots, n-1;& &\|(w_0,w_1)\|_{\mathcal{H}} < \varepsilon. 
  \end{align*}
 We consider the approximated exterior solution 
  \[
   v(x,t) = \sum_{j=1}^n \zeta_j (W_{\lambda_j}(x),0),\qquad (x,t) \in \Omega_0;
  \]
  which satisfies $\|\chi_0 v\|_{Y(\Rm)} \lesssim_1 n$ and solves the following equation in the exterior region $\Omega_0$
  \begin{align*}
   &(\partial_t^2 - \Delta) v = F(v) + e(x,t);& &e(x,t) =  \sum_{j=1}^n \zeta_j F(W_{\lambda_j}) - F(v). 
  \end{align*}
  Here the error term $e(x,t)$ satisfies 
  \begin{align*}
    \|\chi_0 e(x,t)\|_{L^1 L^2} \lesssim_n \kappa. 
  \end{align*}
  If $\kappa$ and $\varepsilon$ are sufficiently small, then an application of perturbation theory (Lemma \ref{perturbation lemma}) implies that the exterior solution $u$ to 
  \[
   \left\{\begin{array}{l} \partial_t^2 u - \Delta u = F(u), \qquad (x,t)\in \Omega_0; \\ (u,u_t)|_{t=0} = (u_0,u_1)\end{array} \right.
  \]
  is defined globally for all $t\in \Rm$ and satisfies 
  \[
   \sup_{t\in \Rm} \|u-v\|_{\mathcal{H}(|t|)} \lesssim_n \kappa + \varepsilon \qquad \Rightarrow \qquad \limsup_{t\rightarrow \pm \infty} \|u\|_{\mathcal{H}(|t|)} \lesssim_n \kappa + \varepsilon. 
  \]
  This implies that the nonlinear radiation profiles $G_\pm$ satisfy $\|G_\pm\|_{L^2(\Rm^+)} \lesssim_n \kappa + \varepsilon$. As a result, we must have $(u_0,u_1) \in \mathcal{R}_n (\delta)$ as long as $\kappa$ and $\varepsilon$ are both sufficiently small. The converse immediately follows from Proposition \ref{main tool}. Here the number of bubbles can be determined by the $\dot{H}^1\times L^2$ norm of $(u_0,u_1)$.
  \end{proof}
 
\begin{remark} \label{enclose explicit}
 A careful review of the proof above shows that given a positive integer $n$, if $\kappa, \varepsilon, \delta$ are sufficiently small, then 
 \begin{align*} 
  &\mathcal{M}_{n}(\kappa, \varepsilon) \subseteq \mathcal{R}_n (c \kappa + c \varepsilon);& & \mathcal{R}_n (\delta) \subseteq \mathcal{M}_n (c \delta, c\delta). 
 \end{align*}
 Here $c$ is a large constant depending on $n$ only. 
\end{remark}
 
 Next we may introduce our ``one-pass theorem''. 

 \begin{proposition} \label{n bubble one pass}
  Assume that $n$ is a positive integer. There exists a positive constant $\delta_0 = \delta_0(n)$, such that if $u$ is a solution to (CP1) with a maximal lifespan $(-T_-,T_+)$ and $\delta$ is a small positive constant $\delta\in (0,\delta_0)$, then the set
  \[
   I_n = \{t\in (-T_-,T_+): \vec{u}(t) \in \mathcal{R}_n (\delta)\}
  \]
  is either empty or an open interval. 
 \end{proposition} 
 \begin{proof}
  We first show that $I = \{t\in (-T_-,T_+): \vec{u}(t) \in \mathcal{R}(\delta)\}$ is either empty or an open interval. Since $\mathcal{R}(\delta)$ is an open subset of $\mathcal{H}$, the continuity of data implies that $I$ is an open subset of $(-T_-,T_+)$. Thus it suffices to show that if $t_1,t_2 \in I$, then  $[t_1,t_2] \subset I$. If $t_1<t_2$ are both contained in $I$, then we may extend the domain of $u$ to 
  \[
   \Rm^3 \times (-T_-,T_+)\cup \{(x,t): |x|>|t-t_1|\} \cup \{(x,t): |x|>|t-t_2|\}.
  \]
  In addition, the time-translated solution $u(x,t+t_1)$ comes with a nonlinear radiation profile $G_{1,-}(s)$ in the negative time direction with $\|G_{1,-}\|_{L^2(\Rm^+)} < \delta$. Similarly the time-translated solution $u(x,t+t_2)$ comes with a nonlinear radiation profile $G_{2,+}(s)$ in the positive time direction with $\|G_{2,+}\|_{L^2(\Rm^+)} < \delta$. It follows that given $t' \in (t_1,t_2)$, the time-translated solution $u(x,t+t')$ is an exterior solution defined in the whole exterior region $\Omega_0$ with initial data $\vec{u}(t')$ and nonlinear radiation profiles 
 \begin{align*}
  & G'_{+}(s) = G_{2,+}(s + t_2 - t'), \; s>0;& & G'_-(s) = G_{1,-} (s+t'-t_1), \; s>0. 
 \end{align*}
 Cleary the inequalities $\|G'_{\pm}\|_{L^2(\Rm^+)} < \delta$ hold, thus $\vec{u}(t') \in I$. In order to finish the proof, we show that either $I_n = \varnothing$ or $I_n = I$ holds, as long as $\delta < \delta_0(n)$ is sufficiently small. Indeed, if $\delta < \delta_0$ is small, Proposition \ref{main tool} implies that 
  \[
  (u_0,u_1) \in \mathcal{R}_n (\delta)\quad \Rightarrow \quad (n-1/3)\|W\|_{\dot{H}^1} < \|(u_0,u_1)\|_{\dot{H}^1\times L^2}^2 < (n+1/3)\|W\|_{\dot{H}^1}. 
 \]
 A continuity argument shows that if $I_n = I \cap I_n \neq \varnothing$, then $I = I_n$. This finishes the proof. 
 \end{proof}
 
 \begin{corollary} \label{n bubble one pass M}
  Given any positive integer $n$ and two constants $\kappa_1, \varepsilon_1 > 0$, there exist two small constants $\kappa_2 < \kappa_1 $ and $\varepsilon_2 < \varepsilon_1$, such that if $u$ is a solution to (CP1) and $t_1,t_2$ are two times satisfying $\vec{u}(t_1), \vec{u}(t_2) \in \mathcal{M}_n(\kappa_2, \varepsilon_2)$, then 
  \[
   \vec{u}(t) \in \mathcal{M}_n (\kappa_1, \varepsilon_1, \{\zeta_j\}), \qquad \forall t \in [t_1,t_2]
  \]
  for some $\{\zeta_j\}_j \in \{+1,-1\}^n$.  
 \end{corollary}
 \begin{proof}
  Without loss of generality, we may assume that $\kappa_1, \varepsilon_1 \ll 1$ such that 
  \[
   \mathcal{M}_n (\kappa_1, \varepsilon_1, \{\zeta_j\}_j) \cap \mathcal{M}_n (\kappa_1, \varepsilon_1, \{\zeta'_j\}_j) = \varnothing, \qquad \{\zeta_j\}_j \neq \{\zeta'_j\}_j. 
  \]
  By Lemma \ref{two equivalent neig}, there exists $\delta = \delta(n, \kappa_1, \varepsilon_1)$, such that $\mathcal{R}_n (\delta) \subseteq \mathcal{M}_n (\kappa_1, \varepsilon_1)$. Without loss of generality we may choose $\delta < \delta_0(n)$. Here $\delta_0(n)$ is the constant given in Proposition \ref{n bubble one pass}. Now we apply Lemma \ref{two equivalent neig} again to find two constants $\kappa_2 < \kappa_1$ and $\varepsilon_2 < \varepsilon_1$ such that $\mathcal{M}_n (\kappa_2, \varepsilon_2) \subseteq \mathcal{R}_n (\delta)$. Now let us assume that $\vec{u}(t_1), \vec{u}(t_2) \in \mathcal{M}_n (\kappa_2, \varepsilon_2)$ and verify that $\vec{u}(t) \in \mathcal{M}_n (\kappa_1, \varepsilon_1, \{\zeta_j\})$ for some $\{\zeta_j\}_j \in \{+1,-1\}^n$. First of all, the inclusion given above implies $\vec{u}(t_1), \vec{u}(t_2) \in \mathcal{R}_n (\delta)$. Thanks to Proposition \ref{n bubble one pass}, we must have 
  \[
   \vec{u}(t) \in  \mathcal{R}_n (\delta) \subseteq  \mathcal{M}_n (\kappa_1, \varepsilon_1), \qquad \forall t\in [t_1,t_2].
  \]
 Finally the existence of $\{\zeta_j\}_j$ follows from the continuity of $\vec{u}(t)$ and the fact that the open sets $\mathcal{M}_n (\kappa_1, \varepsilon_1)$ is a disjoint union of open sets
 \[
  \mathcal{M}_n (\kappa_1, \varepsilon_1) = \bigcup_{\{\zeta_j\}_j \in \{+1,-1\}^n} \mathcal{M}_n (\kappa_1, \varepsilon_1, \{\zeta_j\}_j).
 \]
 This finishes the proof. Finally the author would like to mention that according to Remark \ref{enclose explicit}, if $\kappa_1, \varepsilon_1$ are sufficiently small, then we may choose $\kappa_2 = \varepsilon_2 = c' \min\{\kappa_1, \varepsilon_1\}$, where $c' = c'(n)$ is a small constant.  
 \end{proof}
 
 Before we conclude this section, we characterize all global solutions $u$ to (CP1) defined for all $t\in \Rm$ whose radiation part is small in both two time directions, as an application of our ``one-pass'' theorem given above. 
 
 \begin{corollary}
  Given any $E_0, \kappa,\varepsilon > 0$, there exists a small constant $\delta = \delta(E_0, \kappa, \varepsilon) > 0$ such that if $u$ is a solution to (CP1) satisfying 
  \begin{itemize}
   \item $u$ is defined for all $t\in \Rm$ with an energy $E$ satisfying $E(W,0) \leq E<E_0$; 
   \item The corresponding radiation profiles $G_\pm$, as defined in \eqref{radiation profile u G}, satisfies $\|G_\pm\|_{L^2(\Rm)} < \delta$;
  \end{itemize}
  then there exists a positive integer $J \leq E_0 / E(W,0)$ such that
  \[
   \vec{u}(t) \in \mathcal{M}_J (\kappa, \varepsilon), \qquad \forall t\in \Rm.
  \]
 \end{corollary}
 \begin{proof}
  First of all, we fix a positive integer $n = \lfloor\frac{E_0}{E(W,0)}\rfloor$. We then choose a sufficiently small constant $\delta \ll 1$ such that
  \begin{align*}
   &\delta < \min_{1\leq j\leq n} \delta_0(j); & &\mathcal{R}_j(\delta) \subseteq \mathcal{M}_j (\kappa, \varepsilon), \quad j = 1,2,\cdots, n.
  \end{align*}
  Here $\delta_0(j)$ are the constants in Proposition \ref{n bubble one pass}. Almost orthogonality of the soliton resolution shows that 
  \begin{align*}
   \lim_{t\rightarrow \pm \infty} \|\vec{u}(t)\|_{\dot{H}^1\times L^2}^2 & = 8\pi \|G_\pm\|_{L^2(\Rm)}^2 + J_\pm \|W\|_{\dot{H}^1}^2; \\
   E & = 4\pi \|G_\pm\|_{L^2(\Rm)}^2 + J_\pm E(W,0).  
  \end{align*}
  Here $J_\pm$ are the bubble numbers in the positive/negative time directions. Our assumption on the smallness of $\|G_\pm\|_{L^2(\Rm)}$ implies that $J_+=J_-\in \{1,2,\cdots, n\}$. We let $J=J_+=J_-$. Combining our smallness assumption on the radiation profiles and the limits of $\|\vec{u}(t)\|_{\dot{H}^1\times L^2}^2$ as $t\rightarrow \pm \infty$ given above, we deduce that 
  \[
   \vec{u}(t) \in \mathcal{R}_J (\delta), \qquad |t|\gg 1. 
  \]
  We then apply Proposition \ref{n bubble one pass} and conclude that 
  \[
   \vec{u}(t) \in \mathcal{R}_J (\delta)\subseteq \mathcal{M}_J(\kappa,\varepsilon), \qquad \forall t\in \Rm.
  \]
 \end{proof}
 
\section{Type II blow-up solutions}

In this final section we briefly discuss the type II blow-up solutions. We first sketch a proof of the soliton resolution conjecture for type II blow-up solutions in the 3D radial case, by applying Proposition \ref{main tool}, and then discuss the possibility to evaluate global behaviours of solutions after the blow-up time. 

\subsection{Soliton resolution of type II blow-up solutions} 

Let us consider the soliton resolution conjecture for type II blow-up solutions. Since many intermediate results below have been proved in the previous works, we only give the major steps and skip the details. Let us assume that $u$ is a radial solution to (CP1) defined in an interval $[t_0,0)$ with 
\begin{align*}
 &\|u\|_{Y([t_0,0))} =+\infty;& &M_0 \doteq \sup_{t_0\leq t<0} \|\vec{u}(t)\|_{\dot{H}^1\times L^2} < +\infty. 
\end{align*}
Here $t_0<0$ is a time. Then we may prove the soliton resolution for $t$ close to zero by following the procedure below:
\begin{itemize} 
 \item There exists a finite-energy free wave $v_L = \mathbf{S}_L (v_0,v_1)$ such that 
 \begin{equation*}
  \lim_{t\rightarrow 0^-} \|\vec{u}(t)-\vec{v}_L(t)\|_{\mathcal{H}(|t|)} = 0. 
 \end{equation*}
 In fact, if $v$ is the solution to (CP1) with initial data $(v_0,v_1)$, then the following identity holds if $t<0$ is sufficiently close to zero 
 \begin{equation} \label{canonical evo pre}
  v(x,t) = u(x,t), \qquad |x|>|t|.
 \end{equation} 
 As a result, it suffices to show that given $\kappa, \varepsilon > 0$, there exists a sequence of signs $\zeta_j \in \{+1,-1\}$, and a sequence of scale functions $\lambda_j(t)$ for $j=1,2,\cdots, J$ satisfying 
 \[
  \max\left\{\frac{\lambda_1(t)}{|t|}, \frac{\lambda_2(t)}{\lambda_1(t)}, \cdots, \frac{\lambda_J(t)}{\lambda_{J-1}(t)}\right\} \leq \kappa^2, 
 \]
 such that the following holds for $t<0$ is sufficiently close zero 
 \[
  \int_{|x|<|t|} \left|(\nabla u, u_t) - \sum_{j=1}^J \zeta_j (\nabla W_{\lambda_j(t)}(x),0)\right|^2 {\rm d} x \leq \varepsilon^2. 
 \]
 \item By applying a cut-off technique and using the finite speed of propagation, we may assume that $(v_0,v_1)$ comes with a sufficiently small norm, i.e. $\|(v_0,v_1)\|_{\dot{H}^1\times L^2} < \delta$, where $\delta$ is to be determined later. This means that 
 \[
  \|\vec{v}(t)\|_{\mathcal{H}(|t|)} < 2\delta, \quad t \in \Rm; \qquad \Longrightarrow \qquad \|\vec{u}(t)\|_{\mathcal{H}(|t|)}\leq 2\delta, \quad t\in (t_0,0). 
 \]
 Thus we may fix a small time $t_1<0$ and find a small number $r_1<|t_1|$, such that $\|\vec{u}(t_1)\|_{\mathcal{H}(|t_1|-r_1)} < 3\delta$. As a result, we may extend the domain of $u$ at least to the region 
 \[
  \left(\Rm^3 \times [-t_1,0)\right) \cup \{(x,t): |x|>t>0\} \cup \{(x,t): |x|>|t|-r_1, t<t_1\}. 
 \]
 In addition, the radiation profile of $u$ satisfies 
 \begin{align*}
  &\|G_+\|_{L^2(\Rm^+)} \lesssim_1 \delta;& & \|G_-\|_{L^2([-r_1, +\infty))} \lesssim_1 \delta. 
 \end{align*}
 \item We choose a positive integer $n = n(M_0)$ with $(n-1)\|W\|_{\dot{H}^1}^2 > M_0^2+1$. An application of Proposition \ref{main tool} immediately gives the desired soliton resolution of $\vec{u}(t)$ in the whole space $\Rm^3$ for $t \in (-r_1,0)$ as long as $\delta =\delta(M_0,\kappa,\varepsilon)$ is sufficiently small.  Here we may exclude case (b) by considering the energy norm, as we did in the proof of Proposition \ref{pre soliton resolution}. Please note that  we may prove the inequality 
 \[
  \frac{\alpha_1 (t)^2}{|t|} \lesssim_n \delta^2, 
 \]
 which gives $\lambda_1(t)/|t| \leq \kappa^2$ for small $\delta$, by considering the energy norm outside the main light cone $|x|=|t|$. Finally a continuity argument shows that the bubble number $J$ and the signs of $\alpha_j (t)$ are independent of $t\in (-r_1,0)$. Again we let $\lambda_j(t) = |\alpha_j(t)|^2$ and  $\zeta_j = {\rm sign} \alpha_j(t)$ here. 
\end{itemize}

\subsection{Canonical evolution}

Assume that $u$ is a type II blow-up solution with a blow-up time $T_1$. Let us consider the possible time evolution after the time $T_1$. In view of the identity \eqref{canonical evo pre}, where we apply a time translation if necessary, it is natural to define $u(t)$ to be the solution of (CP1) with initial data $\vec{v}_L(T_1)$ for $t>T_1$. If $u$ blows up again at a time $T_2>T_1$ with 
\[
 \limsup_{t\rightarrow T_2^-} \|\vec{u}(t)\|_{\dot{H}^1\times L^2} < +\infty,
\]
i.e. the blow-up is still of type II, then we may further extend its domain in the same manner. It turns out that this indeed gives a weak solution to (CP1). More details of this type of solutions have been discussed by Krieger-Wong \cite{canonical}. They called this kind of weak solutions ``canonical solutions''. Their results about canonical solutions include the following 
\begin{itemize}
 \item If a canonical solution blows up at a finite time $T_+$, then this blow-up is of type I. Namely 
 \[
  \limsup_{t\rightarrow T_+} \|\vec{u}(t)\|_{\dot{H}^1\times L^2} = +\infty. 
 \]
 \item The functions $J(t)$, $J'(t)$, $J''(t)$ are continuous in the whole lifespan of a canonical solution $u$, although $\|\vec{u}(t)\|_{\dot{H}^1\times L^2}$ is not necessary continuous. Here $J(t)$ is the virtual functional defined in the proof of Proposition \ref{local boundedness}. Please note that Krieger and Wong's  original article defined the function $J(t)$ in a slightly different way from ours, but the same proof still works with minor modifications. 
 \item if the energy $E(t) < 0$ for some time, then the canonical solution $u$ blows up in finite time. Please note that $E(t)$ is not necessary a constant in the whole lifespan of a canonical solution. In fact, if the corresponding strong solution blows up in a manner of type II at some time, then the energy decreases simultaneously by a multiple of $E(W,0)$. 
\end{itemize}
Now let us consider a canonical solution to (CP1). It is not difficult to deduce from the results mentioned above and the virtual identity that unless $u$ blows up at a finite time in the manner of type I, it must satisfy the following
\begin{itemize}
 \item $u$ can be defined for all time $t\geq 0$;
 \item the corresponding strong solution blows up at only finite times $0<T_1<T_2<\cdots<T_m$. The number $m$ can be dominated by the original energy $E$ of $u$
 \[ 
  m \leq \frac{E}{E(W,0)}. 
 \]
\end{itemize}
Let $u$ be such a canonical solution and $R$ be a large radius so that the initial data satisfy $\|(u_0,u_1)\|_{\mathcal{H}(R)} \ll 1$. We may extend the domain of $u$ to the region 
\[
 \left(\Rm^3 \times [0,+\infty)\right) \cup \{(x,t): |x|>R+|t|, t<0\}. 
\]
It is not difficult to see the restriction of $u$ to the exterior region $\{(x,t): |x|>|t-t_1|\}$ is exactly an exterior solution to (CP1), without using the conception of canonical evolution, as long as $t_1 \geq R$. This is true even if $t_1$ is a type II blow-up time for the corresponding strong solution. This is due to the fact that singularity of a canonical solution lies only at the origin and at a type II blow-up time. A careful review of our argument in the previous sections shows that Theorem \ref{main thm} still holds for canonical solutions defined for all time $t\geq 0$ with the following minor modifications: 
\begin{itemize}
 \item Some stable period may be half-closed, half open interval $[a_k,b_k)$, since $b_k$ might be a type II blow-up time of the strong solution. 
 \item Two consecutive stable periods may be separated from each other by either a collision period, with a significant amount of forward radiation, or a single time of type II blow-up (i.e. $b_k = a_{k+1}$). 
\end{itemize}
In other words, the approximated soliton resolution phenomenon given in Theorem \ref{main thm} still occurs for most time $t\geq 0$, except for a few ``short'' intervals, unless the canonical solution eventually blows up in the manner of type I. 

\begin{remark}
 Let $u$ be a canonical solution as above and $T_1$ be a type II blow-up time of its corresponding strong solution. It seems that a few bubbles disappear at the time $T_1$. Accordingly the energy and $\dot{H}^1\times L^2$ norm decreases simultaneously. However, there is another way to understand the scenario at the blow-up time $T_1$. Those bubbles do not disappear, but are permanently ``frozen'' with a fixed scale $\lambda = 0$ after the time $T_1$. The scale separation inequality still holds because 
 \begin{itemize}
  \item If we compare the scales of a frozen bubble with a regular bubble $W_\lambda$, it is clear that $0/\lambda = 0 < \kappa^2$ holds for any $\kappa > 0$; 
  \item If we compare the scales of two frozen bubbles, $0/0$ is an indeterminate form. Thus the inequality $0/0 < \kappa^2$ still holds in a weak sense.
 \end{itemize}
\end{remark}

\section*{Acknowledgement}
The author is financially supported by National Natural Science Foundation of China Project 12471230.

\end{document}